\pdfoutput=1
\documentclass[leqno]{amsart}

\usepackage{amsmath}
\usepackage{amsfonts}
\usepackage{amssymb}
\usepackage{pdfsync}
\usepackage{color,graphicx,tikz-cd}
\usepackage{a4wide}
\usepackage{amscd,amsthm,latexsym}
\usepackage{hyperref}
\usepackage{kotex}
\usepackage{relsize}
\usepackage{cite}

\allowdisplaybreaks %

\usepackage{float}
\usepackage{graphicx}
\usepackage{wrapfig}

\title{A PROOF OF THE FUNCTIONAL EQUATION  CONJECTURE}

\author{Adriano Garsia}
\address{$^1$Department of Mathematics, University of California, San Diego, La Jolla, CA, USA}
\email{garsia@math.ucsd.edu}

\author{Angela Hicks}
\address{ $^2$Department of Mathematics, Lehigh  University,
14 E. Prkwr Ave, Bethlehem, PA 18015 }
\email{anh316@lehigh.edu}

\author{Guoce Xin}
\address{$^3$School of Mathematical Sciences, Capital Normal University, Beijing 100048, PR China}
\email{guoce.xin@163.com}

\date{\today}

\thanks{The first named author was supported by NFS grant DMS1700233.
}



\date{\today}

\theoremstyle{definition}

\def\multi#1{\vbox{\baselineskip=0pt\halign{\hfil$\scriptstyle\vphantom{(_)}##$\hfil\cr#1\crcr}}}
\def \o1  {{\overline }}
\def \l {{\ell}}

\def \BN {{\bf N}}

\def \BC {{\bf C}}

\def \BV {{\bf V}}

\def \bu {\hskip -.15in}

\def \WP {{\widetilde P}}

\def \WV {{\widetilde V}}
\def \WBD  {\widetilde {BD}}

\def\tttt #1{{\textstyle{#1} }}

\def \TAU  {\mathcal{T}}

\def \CL {\mathcal{L}}

\def \uu {\hskip -.08in}

\def \CT {\mathcal{T}}

\def \DA {\downarrow}

\def \magstep#1 {\ifcase#1 1000\or 1200\or 1440\or 1728\or 2074\or 2488\fi\relax}

 \font\Ch=msbm8

\def \UA {\uparrow}

\def\Q{\hbox{\Ch Q}}

\def\la{{\lambda}}

\def \CF {\mathcal{F}}

\font\ita=cmssi10

\font\title=cmbx10 scaled\magstep2

\font\bol=cmbx12

\def \con {\subseteq}

\def\sig{\sigma}

\def \ggg {\gamma}
\def \GG {\Gamma}

\def \-> {\rightarrow}
\def\LL{\big\langle}

\def\RR {\big\rangle}

\def\om {\omega}
\def\la {\lambda}

\def \RA {\rightarrow}

\def \sas {\vskip .06truein}
\def\sa{{\vskip .125truein}}

\def\sapp {{\vskip .5truein}}

\def \eee {\epsilon}

\def\aaa {\alpha}
\def\bbb {\beta}

\def\ggg {\gamma}
\def\aa {\alpha}
\def\bb {\beta}

\def\con {\subseteq}

\def \ses {\enskip = \enskip}
\def \sps {\, + \,}

\def \sms {\, - \,}

\def \scs {\, , \,}
\def \ess {\enskip}

\def \ssp {\hskip .25em}
\def \bigsp {\hskip .5truein}
\def \part {\vdash}

\def \RA {{ \rightarrow }}

\def \CF {\mathcal{F}}

\def \BV {{\bf V}}

\def \BQ {{\Q}}

\def \om {\omega}

\def \TH {{\widetilde H}}

\def \BQ {{\Q}}
\def \om {\omega}

\def \TH {{\widetilde H}}

\def \da {\downarrow}

\def \ua {\uparrow}

\def \scs {\ssp , \ssp}
\def \ess {\enskip}
\def \ssp {\hskip .25em}
\def \bigsp {\hskip .5truein}
\def \part {\vdash}

\def \CL {\mathcal{L}}

\def \TS {\textstyle}

\begin{document}

\vskip -.2in
\begin{abstract}
In the early 2000's the first and second named authors worked for a period of six years in an attempt of proving the Compositional Shuffle Conjecture [1]. Their approach was based on the discovery  that all the Combinatorial properties predicted by the Compositional  Shuffle Conjecture remain valid for each family of Parking Functions with prescribed  diagonal cars. The validity of this property was reduced to the proof of a functional equation satisfied by a Catalan family of univariate polynomials. The main result in this paper is a proof of this functional equation. The Compositional Shuffle Conjecture was proved in 2015 by Eric Carlsson and Anton Mellit [3].
Our proof of the Functional Equation  removes one of the main obstacles in the completion of the Garsia-Hicks approach to the proof of the   Compositional Shuffle Conjecture. At the end of this writing we formulate a few further conjectures including what remains to be proved  to complete this approach.

\end{abstract}


\maketitle

\vspace{-.4 in}
\section{Introduction}

Our manipulations rely heavily on the plethystic notation used in \cite{GXZ}. In fact,  all the notations used in this paper 
is introduced in full detail in
the first section of  \cite{GXZ}.

Recall that Dyck paths in the $n\times n$ lattice square $R_n$ are   paths from $(0,0)$ to $(n,n)$ proceeding by north and east  unit steps, always remaining  weakly above the main diagonal of $R_n$.  These paths are usually represented by their area sequence
$(a_1,a_2,\ldots, a_n)$, where $a_i$ counts the number of complete cells
 between the north step in the $i^{th}$ row and the main diagonal. Notice that the $x$-coordinate of the north step in the $i^{th}$ row is simply the difference $u_i=i-1-a_i$.

 A parking function $PF$ supported by the Dyck path $D\in R_n$ is obtained by labeling the north steps of $D$ with $1,2,\ldots ,n$ (usually referred as ``cars''), where the labels increase along the north segments of $D$. Parking functions can be represented as two line arrays
$$
PF=\begin{pmatrix}
c_1 & c_2  &  \cdots  &   c_n\\
a_1 & a_2  &  \cdots  &   a_n\\
\end{pmatrix}
$$
 with cars $c_i$ and area numbers $a_i$ listed from bottom to top.
We also set
$$
area(PF)= \sum_{i=1}^n a_i,
\ess\ess\ess
dinv(PF)= \uu \sum_{1\le i<j\le n}\uu \Big(
\chi(c_i<c_j \ess\&\ess a_i=a_j)\sps
\chi(c_i>c_j \ess\&\ess a_i=a_j+1)
\Big).
$$
Moreover, the word $w(PF)$ is the permutation obtained by reading the cars
in the two line array by decreasing area numbers and from right to left.

Notice that the diagonals of $R_n$ can be so ordered
that car $c_i$ lies in diagonal $a_i$. Where
diagonal $0$ is the main diagonal of $R_n$.
 A given Dyck path $D$ can hit  the main diagonal in at most $n$ distinct lattice points ($(0,0)$ not counted). We will write $p(D)=p$ for $p= (p_1,p_2,\ldots ,p_{\l(p)})$ a composition  of $n$,
if and only if the  components of $p$ give the sizes of the intervals between successive main diagonal hits of $D$.
If a Parking Function $PF$ is supported by $D$ and $p(D)=p$ it will be convenient to write $p(PF)=p$. This given, we set
$$
\Pi_p[X;q,t]\ses \sum_{p(PF)=p}
 t^{area(PF)}q^{dinv(PF)}F_{ides(w(PF))}[X],
\eqno 1.1
$$
where the last factor in 1.1 is the Gessel Fundamental quasi-symmetric function indexed by the descent set of the inverse of the word $w(PF)$. It is shown in \cite{HHLRU05},  that the right-hand side of 1.1 defines a symmetric function for any composition $p$.
\newpage

If $p= (p_1,p_2,\ldots ,p_{\l(p)})$ then the Compositional Shuffle Conjecture states that the same symmetric function can be obtained
by setting
$$
\nabla C_p[X;q,t]\ses \nabla C_{p_1}C_{p_2}\cdots C_{\l(p)}\, 1
\eqno 1.2
$$
Where $\nabla$ is the eigen-operator of the Modified Macdonald polynomial $\TH_\mu[X;q,t]$ with eigenvalue $T_\mu=t^{n(\mu)}q^{n(\mu')}$. Here, the family
$\big\{\TH_\mu[X;q,t]\big\}_\mu$ is defined as the  unique symmetric function basis that
satisfies the two triangularity  conditions,
(with partition  inequalities indicating dominance)
$$
a)\ess\TH_\mu[X;q,t]=\sum_{\la\le \mu} s_\la \big[\tttt{X\over t-1}\big]c_{\la,\mu}(q,t),
\bigsp\ess\ess
b)\ess \TH_\mu[X;q,t]=\sum_{\la\ge \mu} s_\la \big[\tttt{X\over 1-q}\big]d_{\la,\mu},
\eqno 1.3
$$
together with the normalizing condition
$$
\LL \TH_\mu[X;q,t]\scs s_n\RR= 1.
\eqno 1.3
$$
We can also set for any symmetric function $F[X]$
$$
C_aF[X]\ses \big(-\tttt{1\over q}\big)^{a-1}F\big[X-\tttt{ 1-1/q \over z}\big]\sum_{m\ge 0} z^m h_m[X]\Big|_{z^a}.
\eqno 1.4
$$
It is well known and easily verified from 1.4
that for any pair of positive integers $(a,b)$  we have
$$
q(C_bC_a+C_{a-1}C_{b+1})\ses (C_aC_b+C_{b+1}C_{a-1})
\eqno 1.5
$$
Applying this operator to 1 and then applying  $\nabla$ to the  resulting symmetric function, the equality of the two
polynomials in 1.1 and 1.2, gives the identity
$$
q(\Pi_{b,a}[X;q,t]+\Pi_{a-1,b+1}[X;q,t])
\ses \Pi_{a,b}[X;q,t]+\Pi_{b+1,a-1}[X;q,t]
\eqno 1.6
$$
This identity suggests the existence of a bijection from the family of Parking Functions with diagonal compositions $(a,b)$ or $(b+1,a-1)$ onto the family of Parking Functions with diagonal compositions $(b,a)$ or $(a-1,b+1)$ with the following  properties
\begin{enumerate}
 \item {\it preserves area},
\item {\it preserves  the Gessel Fundamental index},
\item {\it increases dinv exactly by one unit.}
\end{enumerate}

Our proof of the Functional Equation implies the existence of such a bijection even when the family of Parking Functions is restricted by  the  requirement of having a pre-selected collection of cars in each diagonal. This fact suggests that it we may be possible to state and prove the Compositional Shuffle conjecture
as a quasi-symmetric function identity.
This circumstance has been the main driving force in our effort to prove the functional equation.

A {\it legal  schedule} is simply a sequence of integers $W=(w_1,w_2,\ldots,w_{k-1})$ with
 the property that $1\le w_1\le 2$ and $1\le w_i\le  w_{i-1}+1$. The corresponding alphabet is  $X_k=\{x_{-1},x_0,x_1,\ldots ,x_{k-1}\}$.
This given, we construct a family of multivariate polynomials $P_W(X_k;q)$ by the following recursion
\sas

\noindent{\bol Definition 1.1}

{\ita For any legal schedule $W=(w_1,w_2,\ldots ,w_{k-1})$ and any $1\le w\le w_{k-1}+1$ we set
$$
P_{W,w}[X_{k+1};q]\ses
\tttt{x_k-q^w\over 1-q}P_{W}[X_{k};q]
\sps
\tttt{1-x_k\over 1-q}
\Big(P_{W}[X_{k};q]
\Big|_{x_{k-i}\RA qx_{k-i}\, ;\,  1\le i\le w}\ess \Big)
\eqno 1.7
$$
with base case }
 $
P_{\phi}[X_0;q]= qx_{-1}+x_0
 $.
\sas

Our next ingredient is the family of univariate polynomials $Q_W(x;q)$ defined by setting
$$
Q_W(x;q)\ses P_{W}[X_{k};q]\Big|_{x_i=x\, ;-1\le i\le k-1}
\eqno 1.8
$$
\sas

It is easy to  derive from 1.7 that the polynomial $Q_W(x;q)$ has degree $k$ in the $x$ variable and has no constant term, thus it may be written in the form
$$
Q_W(x;q)\ses\sum_{r=1}^k A_r(q) x^r.
\eqno 1.9
$$
This given, computer data led the second named author in \cite{thesis}  to state the following
\newpage

\noindent{\bol Conjecture I}

{\ita For all legal schedules $W=(w_1,w_2,\ldots ,w_{k-1})$,  the coefficients
$A_s(q)$ in the expansion 1.9
satisfy the  identities}
\vskip -.15in
$$
q\big(A_{s+1}(q) +A_{k+1-s}(q)\big)= A_{k-s}(q)+A_s(q)  	
\ess\ess\ess\ess(\hbox{for all $0\le s\le  k $})
\eqno 1.10
$$

An equivalent form of  these identities
may be stated as follows,
\sas

\noindent{\bol Conjecture II}

{\ita For all legal schedules $W=(w_1,w_2,\ldots ,w_{k-1})$ the polynomials $Q_W(x;q)$ satisfy the following
functional equation}\vskip -.3in
$$
(1-\tttt{q\over x})Q_W(x;q)+x^k(1-qx)
Q_W(1/x;q)= (1+x^k) (A_k -q A_1 )=(1+x^k)(1-q^2)\prod_{i=1}^{k-1}[w_i]_q
\eqno 1.11
$$
Since 1.11 is invariant under the replacement $x\RA 1/x$, it follows that, to  prove this equivalence, we only need to show that the equalities in 1.10 are obtained  by equating coefficients of $x^s$ on both sides of 1.11 for $0\le s \le k$.

To this end, notice that the coefficient of $x^s$ in the left  hand side of 1.11 is
$$
A_s-q A_{s+1}\sps A_{k-s} -q A_{k+1-s}.
$$
and this must vanish for $1\le s\le k-1$ if 1.11 must hold true.
Now for $s=0$ this reduces to $-qA_1+A_k$, which is precisely what the right hand side of 1.11 gives. Finally if $s=k$ we obtain
$A_k-A_1$ which is again what the right hand side of 1.11 gives.
The last equality in 1.1 follows by setting $x=1$ in the left hand side of 1.11 and using formula 1.19 proved later.
\sas
\noindent{\bol Remark 1.1}

It should be apparent that 1.6 and 1.10 are closely related. In the last section this fact will play a crucial role in conveying the significance of our proof of the functional equation.
\sas

Our next task is the introduction of the basic tool that will be used in our proof. These are the so called ``Bar Diagrams''. This tool was created in \cite{thesis} precisely for this purpose. In fact, the device has so far been very effective   proving special cases of the functional
equation. Roughly speaking  Bar Diagrams  are none other than combinatorial structures that  give  a visual representation of the terms of the polynomials $P_W(X_k;q)$.
More precisely, each of these polynomials can be written in the form
$$
P_W(X_k;q)\ses \sum_{S\con X_k}A_S(q)\,
\prod_{x_i\in S}x_i
\eqno 1.12
$$
\vskip -.03in
\noindent
where the sum is over subsets of the alphabet
$X_k=\{x_{-1},x_0,x_1,\ldots ,x_{k-1})$ with the following  special properties. Given the schedule $W=(w_1,w_2,\ldots ,w_{k-1})$,
define for  $1\le i\le k-1$
$$
act(x_i)\ses \big\{x_{i-1},\ldots,x_{i-w_i}\big\}.
\eqno 1.13
$$
This given, we have, for each $S$ in 1.12
$$
(i)\ess  \hbox{ {\it For each $x_i\in S$, $|S\cap act(x_i)|\ge 1$}},
\ess\ess\ess\ess\ess\ess\ess\ess
(ii)\ess \hbox{\it {\it For each $x_i\in S^c$, $| S^c\cap act(x_i)|\ge 1$ }}
$$
$$
(iii)\ess\ess\hbox{\it $|\{x_{-1},x_0\}\cap S|=1$.}
$$
where $S^c$ denotes the complement of $S$ in $X_k$.

To better understand the mechanism that yields
the coefficient $A_S(q)$ in 1.12  we need
to introduce the notion of ``Labelled Bar Diagram'' corresponding to a given schedule $W$.
By summing over all these diagrams we will obtain a polynomial $\WP_W(Y_k;q)$ in the non-commutative alphabet
$Y_k=\big\{y_{-1},y_0,y_1,\ldots,y_{k-1}\big\}$
with expansion
$$
\WP_W(Y_k;q)\ses \sum_{LBD\in {\mathcal LBD}_W}
m_{LBD}(q)\,\om_{LBD}(Y_k).
\eqno 1.14
$$
Here the sum is  over the family of $LBD's $ constructed  from $W$,   $\om_{LBD}(Y_k)$ is an injective word in the alphabet $Y_k$ and $m_{LBD}(q)$
is  a monomial in $q$ depending only on $LBD$. Moreover, if we denote by  $\om_{LBD}(X_k)$ the result of replacing each letter $y_i$ in $\om_{LBD}(Y_k)$ by the corresponding letter $x_i$ then there is a subset  $S\con X_k$ satisfying properties (i) and (ii)   such that
$\om_{LBD}(X_k)=\prod_{x_i\in S}x_i$.
\newpage

In particular it will follow  that the coefficient $A_S(q)$ in 1.12
can be expressed as
$$
A_S(q)\ses
\sum_{\om_{LBD}(X_k)=S }m_{LBD}(q).
\eqno 1.15
$$

\noindent
A Labelled Bar Diagram of the polynomial $\WP_{2,3,4,4,3,4}(Y_k;q)$ is constructed as depicted  below
\vskip -.15in
\begin{figure}[H]
\centering
\includegraphics[width=4.3in]{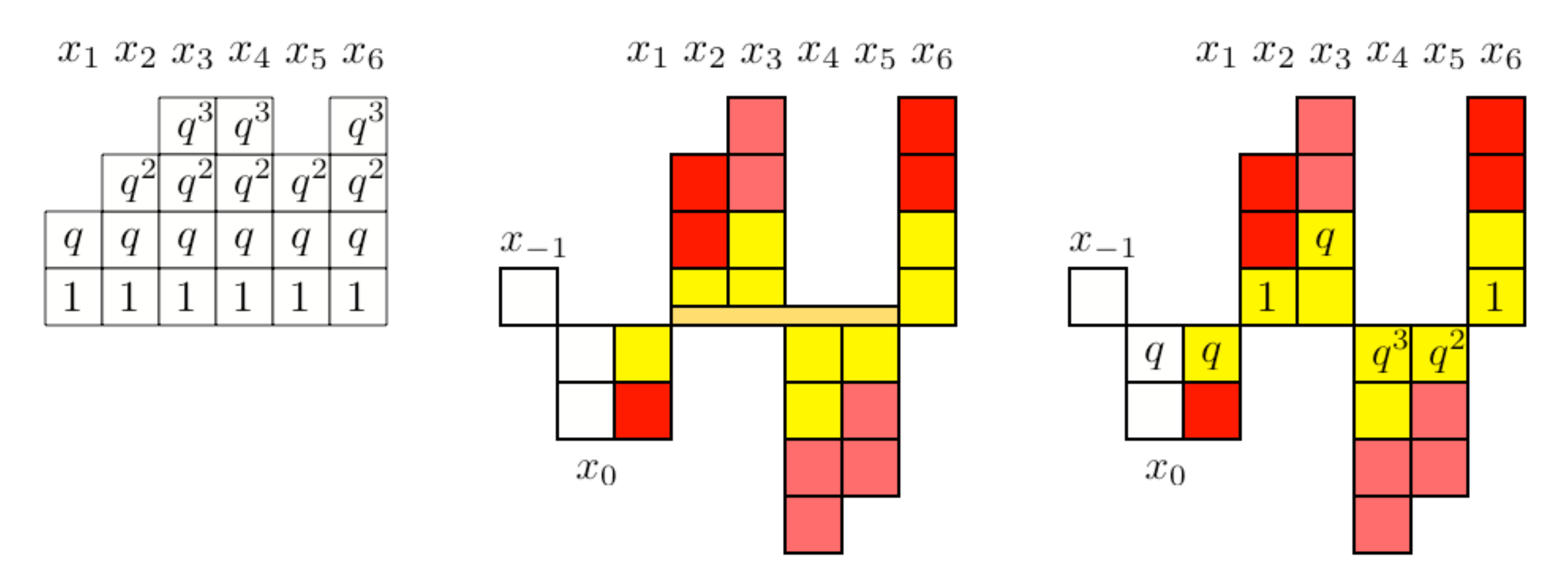}
\end{figure}

\vskip -.25in
\noindent
In the left of this display we have a table of powers of $q$, for each bar. This table will be used for labeling each bar, whether or not the bar is placed above or below the ground line. In the middle of the display we illustrate the mechanism that is used to color the cells of each bar. If we rotate counterclockwise by $90^o$   bar   $x_6$, we see that it will touch two  previous bars. Here we have depicted by the thin rectangle the final position of the rotated bar.
Accordingly,  the first two cells are colored yellow and the cells above them are colored red. All the other bars above the ground line
are colored by the same mechanism. For the bars below the ground line the mechanism is the same, except that in this case we rotate the bar by $90^o$ clockwise. As a result we see that bar $x_4$, after  this rotation, will touch two previous bars. Accordingly, the  two cells
closest to the ground line are colored yellow and the remaining lower cells are colored red. The  subset $S\con X_k$ corresponding to the resulting Bar Diagram is obtained by placing $x_i$ in $S$ if and only if its bar is above the ground line. This given, the only additional property we will require is that in all
our Bar Diagrams each bar has at least one yellow cell.
This is to   assure the   properties (i) and (ii).
The way we draw bar  $x_{-1}$ and bar $x_0$ assures that only  one of them will be up in all Bar Diagrams.  This to assure
property (iii).

It remains to describe how the labeling is done. Notice first  that in bars $x_1,x_2,x_5$ there is only one yellow cell. In a labelled Bar Diagram,
each bar must be given a label. The labels of the yellow cells of these bars, are obtained from the table. Next notice
that, for our example, in each of bars $x_3,x_4,x_6$ there are
two yellow cells. Since we are only allowed one label per bar we have $2\times 2\times 2  $ possible choices here. We displayed only one of them. However, for each of the bars, the chosen cell must be  labeled according to the table.
Finally, notice the $q$ in the top cell of bar $x_0$. This will be always the case whenever
bar $x_{-1}$ is up and bar $x_0$ is down and only then.

\begin{wrapfigure}[8]{R}{1.3in}
\vspace{-14pt}
\includegraphics[width=1.2in]{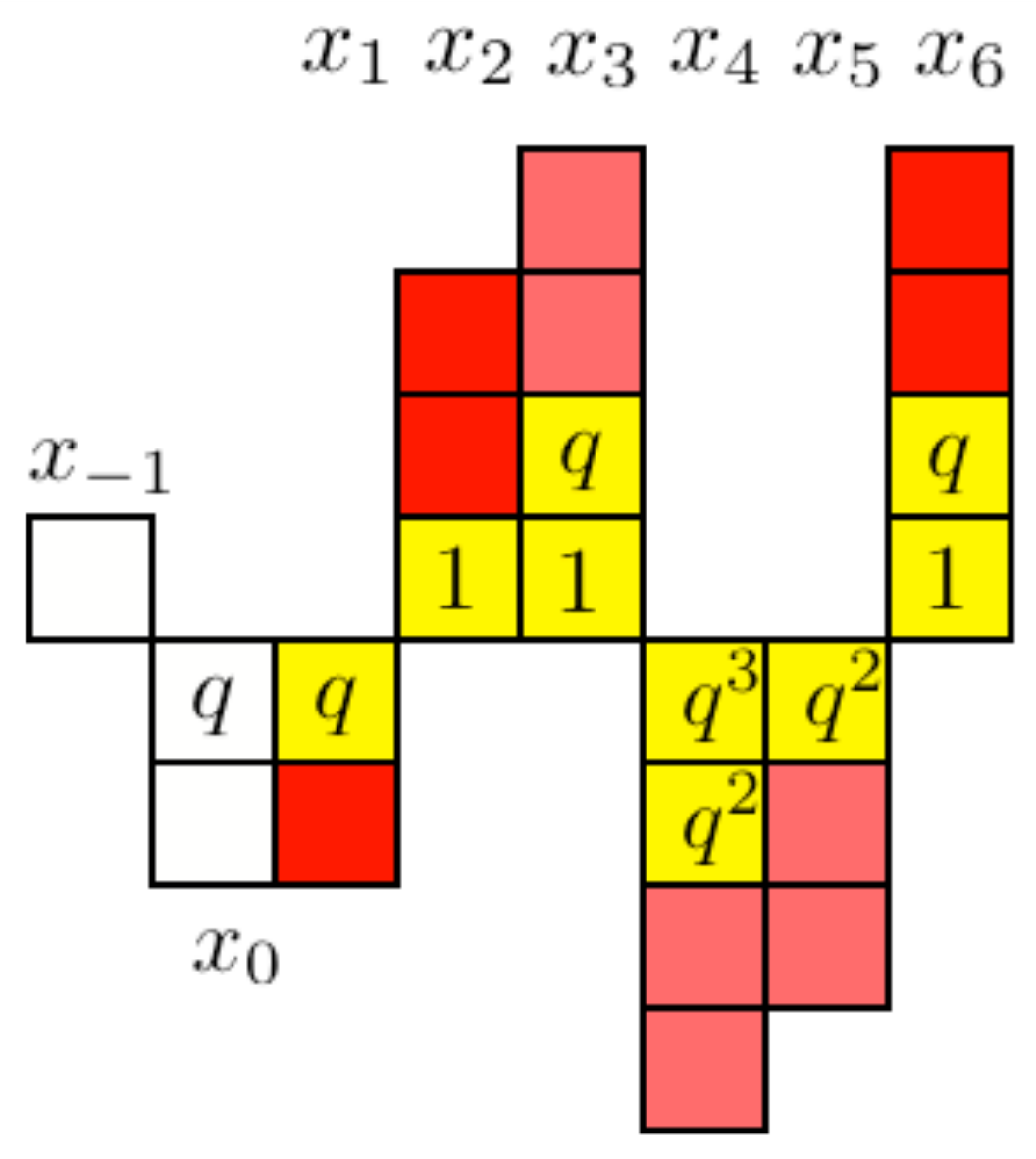}
\end{wrapfigure}

Next, we  need to describe how the word $\om$ is  constructed. To begin, the word $\om$ must  be constructed one letter at the time starting with $y_{-1}$ or $y_0$ according as  which of bars $x_{-1}$ or bar $x_0$ is up. Next, a letter  $y_i$ appears in $\om$ if and only if bar $x_i$ is up. This given, the  power of $q$ used to label  bar $x_i$ dictates where in $\om$ the letter $y_i$ is to be inserted. For instance, if the label is $q^k$ then $y_i$ must be inserted exactly to the left of $k$ from the letters $y_{i-1},y_{i-2},\ldots, y_{i-w_i}$.
The resulting monomial  $m_\om(q)$ is the product  of the powers of $q$ used in the labels.

We will see in the final section of this paper that Labelled Bar Diagrams  are in bijection with Parking Functions with composition of length $2$. Our goal here is to use un-labeled Bar Diagrams constructed for the schedule $W=(w_1,w_2,\ldots, w_{k-1})$
to represent the terms in the expansion
of the polynomial $P_W(X_k;q)$. The idea can be communicated by a single example. We simply  label all the empty yellow cells with the power of $q$ as dictated by the table. As illustrated in the above display. This given, the term of $P_{2,3,4,4,3,4}(X_k;q)$ produced by this un-labeled Bar Diagram $BD$   is none other than
$$
 A_{BD}(q)\, m_{BD}(X_k)\ses q^4(1+q)(q^2+q^3)(1+q)
 \,\, x_{-1}x_2x_3x_6.
$$

The following result shows that we can use labelled bar diagrams to geometrically represent our polynomials $P_W(X_k;q)$ and ultimately prove results about the polynomials $Q_W(x;q)$.
\sas

\noindent{\bol Theorem 1.1}

{\ita For all legal schedules $W=(w_1,w_2,\ldots,w_{k-1})$ we have

\vskip-.14in$$
P_W(X_k;q)\ses \sum_{ BD\in {\mathcal  BD}_W}
  A_{BD}(q)  \,\,m_{ BD}(X_k)
\eqno 1.16
$$
\vskip-.05in
\noindent
where ${\mathcal BD}_W$ denotes the family of all  bar diagrams corresponding to the schedule $W$,
with  $ A_{BD}(q)$ and $m_{ BD}(X_k)$ represent the polynomial in $q$ obtained from  $BD$, and the  monomial   $\prod_{ x_i\in S} x_i$ obtained by letting $S$ be the subset of $X_{k-1}$ of the elements whose bars are up in $BD$}

\noindent{\bol Proof}

We need only show that the polynomial on the right hand side of 1.16 can be obtained by the recursive algorithm of definition 1.1. Since   the base case 1.16 reduces to
$$
P_{\phi}[X_0;q]= q\, x_{-1}+x_0,
$$
we will proceed by induction on the length of the schedule. Let us assume that 1.16 is valid for schedules of length $k-1\ge 0$. Now given $W=(w_1,w_2,\ldots,w_{k-1})$
and any integer $1\le w\le w_{k-1}+1$, we can obtain all the un-labelled bar diagrams in $\mathcal BD_{W,w}$ by starting with an un-labeled  Bar Diagram in $BD\in \mathcal BD_{W}$ and appending to it first an $x_k$ up bar of length $w$ and then an $x_k$ down bar of same length. Calling the resulting Bar diagrams $\WBD^{(1)}$ and $\WBD^{(2)}$. To determine the contributions that these two Bar Diagrams yield to the polynomial
$$
\sum_{\WBD\in {\mathcal  BD_{W,w} }} A_{\WBD}(q)  \,\,
m_{ \WBD}(X_{k+1})
$$
we only need to know one integer. Namely, if
$m_{BD}(X_k) =\prod_{x_i\in S}x_i$ then what  we need is the cardinality
\vskip -.2in
$$
a=\big|S\cap \big\{x_{k-1},x_{k-2},\ldots ,x_{k-w}\big\}\big|
\eqno 1.17
$$
This given,  we see that an up bar of length $w$ appended at the end of the Bar Diagram $BD$ will necessarily have exactly $a$ yellow cells. For the same reason a down bar of length $w$ appended at the end of $BD$ will have exactly $w-a$ yellow cells. This gives, for $0< a< w$
\begin{align*}
A_{BD^{(1)}}(q)\, m_{BD^{(1)}}(X_{k+1})
&= A_{BD}(q) \, m_{BD}(X_k)\,\,
x_k(1+\cdots +q^{a-1})
=  A_{BD}(q) \, m_{BD}(X_k)\,\,\tttt{x_k-x_k q^a\over 1-q}
\\
A_{BD^{(2)}}(q)\, m_{BD^{(2)}}(X_{k+1})
&=
A_{BD}(q) \, m_{BD}(X_k)\,\,(q^a+ \cdots+ q^{w-1})=
A_{BD}(q) \, m_{BD}(X_k)\,\,\tttt{q^a-q^w\over 1-q}.
\end{align*}
In case $a=w$ or $a=0$ then one of the two terms vanishes, but we can still write
\begin{align*}
A_{BD^{(1)}}(q)\,& m_{BD^{(1)}}(X_{k+1})
+
A_{BD^{(2)}}(q)\, m_{BD^{(2)}}(X_{k+1})
\\
&\ses
A_{BD}(q) \, m_{BD}(X_k)
\big(
\tttt{x_k-x_k q^a\over 1-q}
+\tttt{q^a-q^w\over 1-q}
\big)
=
A_{BD}(q) \, m_{BD}(X_k)
\big(
\tttt{x_k- q^w\over 1-q}
+q^a\tttt{1-x_k\over 1-q}
\big)
\\
&\ses
\tttt{x_k- q^w\over 1-q}
\,\,
A_{BD}(q) \, m_{BD}(X_k)
\sps
\tttt{1-x_k\over 1-q}
A_{BD}(q) \, m_{BD}(X_k)
\Big|_{x_{k-i}\RA qx_{k-i}\, ;\,  1\le i\le w}
\end{align*}
It follows from this that by summing over all $BD\in {\mathcal BD}_W $ we will necessarily obtain the identity
\vskip -.2in
\begin{align*}
\sum_{\WBD\in {\mathcal  BD_{W,w} }} \bu A_{\WBD}(q) & \,\,
m_{ \WBD}(X_{k+1})
\ses
\\
\ses &
\tttt{x_k- q^w\over 1-q}\bu\sum_{BD\in {\mathcal  BD_{W} }}\bu A_{BD}(q)  \,\,
m_{ BD}(X_{k})
\sps
\tttt{1-x_k\over 1-q}\sum_{BD\in {\mathcal  BD_{W} }}\bu A_{BD}(q)  \,\,
m_{ BD}(X_{k})\Big|_{x_{k-i}\RA qx_{k-i}\, ;\,  1\le i\le w}
\end{align*}
By the inductive hypothesis, this is none other than
$$
\sum_{\WBD\in {\mathcal  BD_{W,w} }} \bu A_{\WBD}(q)  \,\,
m_{ \WBD}(X_{k+1})
\ses
\tttt{x_k- q^w\over 1-q}
P_W(X_k;q)
\sps
\tttt{1-x_k\over 1-q}
P_W(X_k;q)
\Big|_{x_{k-i}\RA qx_{k-i}\, ;\,  1\le i\le w}
\eqno 1.18
$$
and  Definition 1.1 gives
$$
\sum_{\WBD\in {\mathcal  BD_{W,w} }} \bu A_{\WBD}(q)  \,\,
m_{ \WBD}(X_{k+1})\ses P_{W,w}[X_{k+1};q],
$$
completing the induction and the proof.
\newpage
\noindent{\bol Remark 1.2}

An immediate by-product of this proof is the following identity
$$
P_W[X_k;q]\Big|_{x_i=1;\forall\, i}
\ses (1+q)\prod_{i=1}^{k-1}[w_i]_q
\eqno 1.19
$$
In fact it follows from 1.18, by setting $x_k=1$, that we have
$$
P_{W,w}[X_k;q]\Big|_{x_i=1;\forall\, i}\ses [w]_qP_W[X_k;q]\Big|_{x_i=1;\forall\, i}
$$
This is precisely the  step need to prove 1.19  by induction starting from the base case.
\sa

The following result gives an entirely explicit expression for the polynomials  $A_S(q)$ in 1.12.
\sas

\noindent{\bol Theorem 1.2}

{\ita For a  schedule $W=(w_1,w_2,\ldots ,w_{k-1})$ and a given
$BD\in {\mathcal BD}_W$ let $m_{BD}=\prod_{x_i\in S}x_i$ and let
$a_i=a_i(BD)$ denote the number of $x_{i-1},x_{i-2},\ldots, x_{i-w_i}$ that are in $S $ then }
$$
A_S(q)\, =\,  q^{\chi(x_{-1}\in S )}
\prod_{x_i\in S }[a_i]_q
\times
\prod_{x_i\not\in  S }q^{a_i}[w_i-a_i]_q
\eqno 1.20
$$
\noindent{\bol Proof}

This is one of the by-products of the proof of
Theorem 1.1. See the equalities after the display in 1.17.
\sa

This paper has three more sections. In the next section we gather all  the auxiliary results we need for our proof. In the third section we give the proof of the Functional Equation. In the last section we state what remains to be proved and state some further conjectures.
\sa

\section{Auxiliary results}

This section contains the minimal set of facts we   need from \cite{thesis} for the  proof of the Functional Equation. We will include their  proofs for sake of completeness
\sas

\centerline {\ita The Complementation identity}
\sas
A close look at a single labelling $LBD$
of a bar diagram $BD$ that contributes to
the sum in 1.14 reveals that
the polynomial $P_W(X_k;q)$ possesses
  a degree flipping involution onto itself. This useful fact will play a significant role in establishing various properties of these polynomials. In the display below we have depicted an $LBD$   for the schedule $
W = (2, 3, 2, 3, 3, 4, 3, 4, 5, 5, 3, 4)
$ together with  the labelled diagram obtained by
flipping  LBD  across the ground line
\vskip -.1in
\begin{figure}[H]
\centering
\includegraphics[width=4.5in]{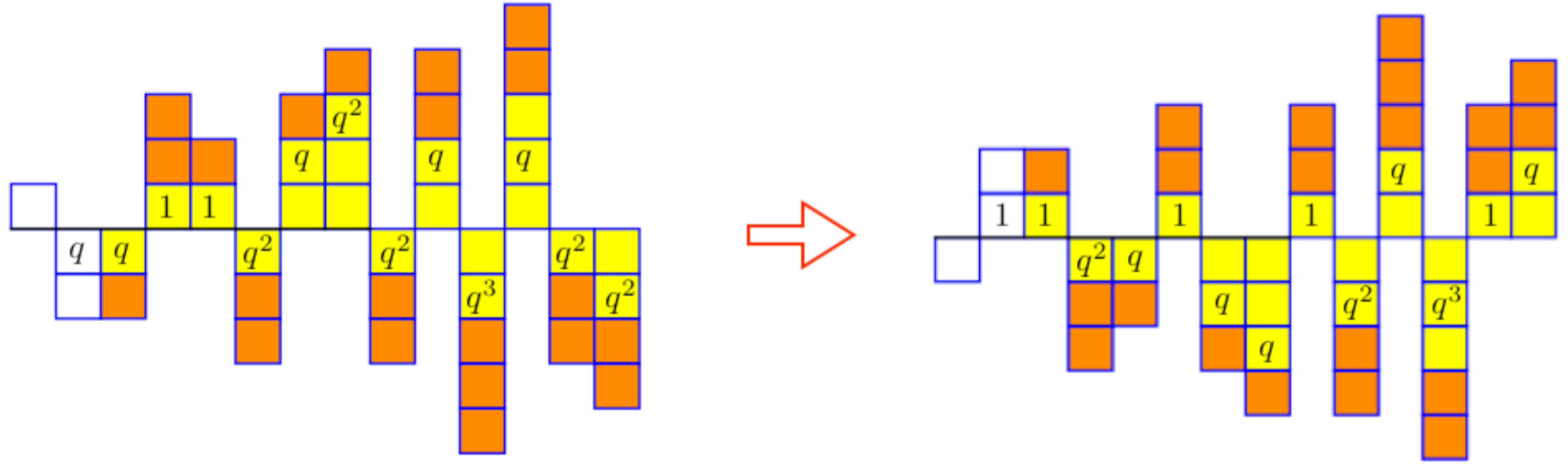}
\end{figure}

\vskip -.2in
We can immediately see, in the above display,  that if for bar $x_i$
of LBD the label is $q^{r_i}$ then the label of bar $x_i$ in $flip(LBD)$ is necessarily  $q^{s_i}$ with  $r_i+s_i=w_i-1$.  This simple example makes it evident that we have a involution of the LBD's which complements the monomials
 $m_{BD}[X_k]=x_\eee\prod_{j=1}^{k-1}x_j^{\chi(bar\,\,j\,\, is\,\, up)}$ , as subsets of the alphabet $X_k$ and complements  the power of $q$ giving the total  the weight of $LBD$.
More precisely we get for some $\eee \in \{-1,0\}$
$$
weight_{flip(BD)}[X_k]=
{x_{-1}\cdots x_{k-1}\over x_\eee}
\prod_{j=1}^{k-1}{1\over x_j^{\chi(bar\,\,j\,\, is\,\, up)}}\, \times\, {q^{1+\sum_{i=1}^{k-1}(w_i-1)}
\over q^{\chi(x_\eee=x_{-1})}} \prod_{i=1}^{k-1}{1\over q^{r_i(LBD)}}
$$
 When these identities  are summed over all $LBD$ we obtain a simple proof of the following  basic result

\newpage

\sas

\noindent{\bol Theorem 2.1} {\rm (The  complementation identity)} 

{\ita For all legal schedules $W=(w_1,w_2,\ldots ,w_{k-1})$ and
$X_k=\{x_{-1},x_0,\ldots ,x_{k-1}\}$}
$$
 P_W(X_k;q)
\ses
x_{-1} \cdots x_{k-1}q^{1+\sum_{i=1}^{k-1}(w_i-1)}
P_W(\tttt{{1\over x_{-1}},{1\over x_0},{1\over x_1},\ldots{1\over x_{k-1}};{1\over q}})
\eqno 2.1
$$

As   immediate corollary we have
\sas

\noindent{\bol Theorem 2.2}

{\ita For all legal schedules $W=(w_1,w_2,\ldots ,w_{k-1})$
the coefficients in the expansion
$$
Q_W(x;q)=\sum_{b=1}^k x^b A_b(q)
$$

\vskip -.25in
\noindent
satisfy the equalities}
$$
A_b(q)\ses q^{1+\sum_{i=1}^{k-1}(w_i-1)}
A_{k+1-b}(1/q)
\bigsp\hbox{(for $1\le b\le k$)}
\eqno 2.2
$$
\noindent{\bol Proof}

Since by definition
$$
Q_W(x;q)=P_W[X_k;q]\Big|_{x_i=x; \forall i}
$$
from 2.1 it follows that
$$
Q_W(x;q)\ses q^{1+\sum_{i=1}^{k-1}(w_i-1)}
x^{k+1}Q_W(1/x;1/q).
$$
Equating coefficients of $x_b$ gives 2.2.
\sa

\centerline {\ita The validity of the Functional equation when adding a component $1$ to the schedule}
\sa

The power of bar diagrams can be gauged
by the following
surprising general fact.
\sas

\noindent{\bol Theorem 2.3}

{\ita  For any legal schedule $W=(1,w_2,\ldots,w_{k-1})$ the polynomial
$Q_W(x;q)$ satisfies  the functional equation. In particular we have}
$$
Q_W(x;q)= A_1(q)x+A_k(q)x^k
\bigsp\big(\hbox{with $A_1(q)=q A_k(q)$}\big)
\eqno 2.3
$$
\noindent{\bol Proof}

In the display below we have depicted 6 diagrams when   $w_1=1$. On the left we have examples of the case when bar $x_{-1}$ is up and on the right we have the case when bar $x_0$ is up.

\vskip -.1in
\begin{figure}[H]
\centering
\includegraphics[width=5.5in]{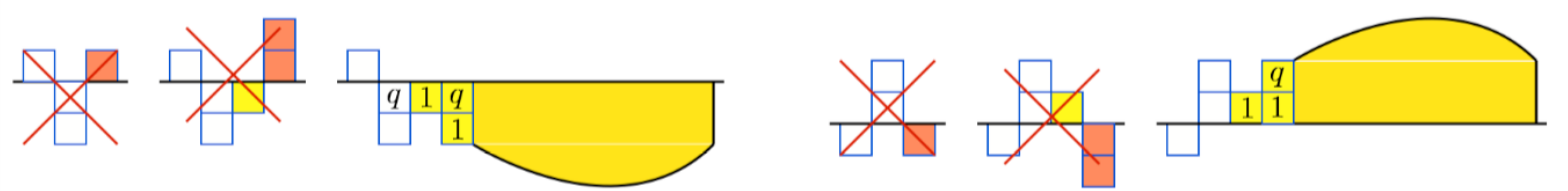}
\end{figure}
Notice that since $x_1$ only acts on $x_0$ we see, in both cases, that if bar $x_1$ is not on the same side as bar $x_0$ we end up with an illegal coloring. But even if bar $x_1$ is on the same side but bar $x_2$ is not then again we end up with an illegal coloring,
(remember that $w_i\le w_{i-1}+1$). In conclusion we see that in any cases all
bars $x_1,\ldots,x_{k-1}$ must be on the same side as bar $x_0$, This proves 2.3.

Let us now substitute 2.3 into the functional equation. This gives
$$
A_1x+A_kx^k\sms qA_1-qA_k x^{k-1}
\sps A_1x^{k-1}+A_k\sms q A_1x^k -qA_k x\ses (1+x^k)(A_k-qA_1)
$$
Thus we see that the functional equation will be  satisfied if and only if $A_1=qA_k$ and $A_k=\prod_{i=1}^{k-1}[w_i]_q$. But one look at the labeling in the above display will clearly show that these two conditions are trivially always satisfied.
\sas

Theorem 2.3 shows there is no loss restricting all our legal schedules to start with $w_1=2$. This is what we will do for the rest of this paper. The next result, proved in \cite{thesis}, shows that the same conclusion  can be drawn even when a later component happens to be equal to $1$.
\newpage

\noindent{\bol Definition} 

{\ita Let us say that a legal schedule $W$ is tame if and only if $Q_{W }(x;q)$ satisfies the functional equation.}

Our next aim here is to show that if two legal schedules  $W'=(w_1,\ldots,w_{k-2}$ and
$W =(w_1,\ldots,w_{k-2},w_{k-1})$ are both tame
then all the legal schedules
$(w_1,\ldots,w_{k-1},1,w_{k+1},\ldots,w_{\l-1})$ will be tame.
It will, be convenient here to introduce the operator $B_{k,w}$ whose action on a polynomial
$A(X_k,q)$ is defined by setting, for $1\le w\le w_{k-1}+1$
$$
B_{k,w}A(X_k,q)\ses \tttt{x_k -q^w\over 1-q}A(X_k,q)\sps \tttt{1 -x_k\over 1-q}
A(X_k,q)\Big|_{x_{k-i}=qx_{k-i}\, ;\,  1\le i\le  w}.
\eqno 2.15
$$
Next, given a legal schedule $W=(w_1,w_2,\ldots,w_{k-1})$ we have  the natural decomposition
$$
P_W(X_k;q)\ses A_W(X_{k-1},q)+x_{k-1}B_W(X_{k-1},q).
\eqno 2.16
$$
This induces the decomposition
$$
Q_W(x;q)\ses A_W(x ,q)+x B_W(x;q),
\eqno 2.17
$$
where by a slight abuse of notation we have set
$$
a)\ess A_W(x ,q)= A_W(X_{k-1})\Big|_{x_i=x;\forall i},
\bigsp
b)\ess B_W(x ,q)= B_W(X_{k-1})\Big|_{x_i=x;\forall i}.
\eqno 2.18
$$
It is clear (from 2.17) that both
$A_W(x;q)$ and $B_W(x;q)$ are of $x$-degree at most $k-1$, but we can be more precise.
\sas

\noindent{\bol Proposition 2.1}

{\ita Setting $d=1+\sum_{i=1}^{k-1}(w_i-1)$ we have
\begin{align*}
&a)\ess\ess\ess q^dx^kB_W(1/x;1/q)\ses A_W(x;q),
\\
&b)\ess\ess\ess q^dx^kA_W(1/x;1/q)\ses B_W(x;q).
\\
\end{align*}

\vskip -.54in
\hfill $2.19$

\vskip  .08in
\noindent
From this it follows that
$$
(1)\ess  B_W(x;q) \hbox{is of x-degree exactly $k-1$},
\bigsp\bigsp
(2)\ess
   A_W(x;q)\big|_x\neq 0
$$
}
\noindent{\bol Proof}

From the complementation identity it follows that
\vskip -.2in
\begin{align*}
q^d\prod_{i=-1}^{k-1}x_i\Big(
A_W(\tttt{1\over x_{-1}} ,\ldots, \tttt{1\over x_{k-2}};\tttt{1\over q})
+ \tttt{1\over x_{k-1}} B_W(\tttt{1\over x_{-1}} ,\ldots, \tttt{1\over x_{k-2}};\tttt{1\over q})\Big)
 =
A_W(X_{k-1};q)+ x_{k-1}B_W(X_{k-1};q).
\end{align*}
Expanding both sides and setting $x_{k-1}=0$ gives
$$
q^dx_{-1}\cdots x_{k-2}B_W(x_{-1}^{-1},\ldots x_{k-2}^{-1};1/q)
\ses A_W(X_{k-1};q)
\eqno 2.20
$$
while equating coefficients of $x_{k-1}$ we get
$$
q^dx_{-1}\cdots x_{k-2}
A_W(x_{-1}^{-1},\ldots x_{k-2}^{-1};1/q)
\ses
B_W(X_{k-1};q)
\eqno  2.21
$$
Thus 2.19 a) and b) follow by replacing all the $x_i$ by $x$ in 2.20 and 2.21.
Note next that since the insertion of
$x_{k-1}$ is what causes the x-degree of $Q_W(x;q)$ to reach $k$,
then  $x$-degree of $B_W(x;q)$ must be exactly $k-1$. This given, 2.19 a)
yields that we must also have $A_W(x;q)\big|_x\neq 0$.
\sas

We are now ready to present  a first surprising fact.

\noindent{\bol Theorem 2.4}

{\ita Let $W=(w_1,\ldots, w_{k-1})$ be a schedule and  $1\le w\le w_{k-1}+1$, then
 both  $Q_W(x;q)$  and $Q_{W,w}(x;q)$ will satisfy the functional equation if and only if

\vskip -.2in
$$
(1-q/x)  A_{W,w}(x;q)
\sps x^{k}(1-xq)B_{W,w}(1/x;q)  \ses
 [w]_q(1-q^2)\prod_{i=1}^{k-1}[w_i]_q
\eqno 2.22
$$

\vskip -.08in
\noindent
In particular, by the complementation result,  we must also have}

\vskip -.18in
$$
(1-q/x)  B_{W,w}(x;q)
\sps x^{k}(1-xq)A_{W,w}(1/x;q)  \ses
x^k [w]_q(1-q^2)\prod_{i=1}^{k-1}[w_i]_q
\eqno 2.23
$$
\noindent{\bol Proof }
\newpage

\noindent
The idea is to start by writing
$$
P_W(X_k;q)\ses \sum_{r=0}^w A_r(X_k;q)
\eqno 2.24
$$
where $A_r(X_k;q)$ is the sum  of all terms in $P_W(X_k;q)$ that contain exactly $r$ of the variables acted upon  by $x_k$.
Setting all $x_i=x$, by a slight abuse of notation we will also write
$$
Q_W(x;q)=\sum_{r=0}^wA_r(x;q)
\eqno 2.25
$$

\vskip -.15in
\noindent
with
$$
A_r(x;q)\ses A_r(X_k;q)\Big|_{x_i=x\, ; \forall i }
$$

\vskip -.1in
$$
P_{W,w}(X_k;q)\ses
\sum_{r=0}^w A_r(X_k;q)
\big( [r]_q x_k +  q^{r} [w-r]_q\big)
\eqno  2.26
$$

\vskip -.15in
\noindent
thus we may write
$$
Q_{W,w}(x;q)\ses
\sum_{r=0}^w A_r(x;q)
\big( [r]_q x \sps  q^{r} [w-r]_q\big)
\eqno 2.27
$$
Notice,  that 2.26 now yields the decomposition

\vskip -.1in
$$
P_{W,w}(X_k;q)\ses
A_{W,w}(X_k;q)+ x_k B_{W,w}(X_k;q)
\eqno 2.28
$$
with
$$
A_{W,w}(X_k;q)= \sum_{r=1}^w A_r(X_k;q)q^{r} [w-r]_q
\ess\ess\ess\ess\ess\hbox{and}
\ess\ess\ess\ess\ess
B_{W,w}(X_k;q)= \sum_{r=1}^w A_r(X_k;q) [r]_q
\eqno 2.29
$$
and

\vskip -.3in
$$
Q_{W,w}(x;q)\ses
A_{W,w}(x;q)\sps xB_{W,w}(x;q)
\eqno 2.30
$$
with

\vskip -.2in
$$
A_{W,w}(x;q)=  \sum_{r=1}^w A_r(x;q)q^{r} [w-r]_q
\ess\ess\ess\ess\ess\hbox{and}
\ess\ess\ess\ess\ess
B_{W,w}(x;q)=\sum_{r=1}^w A_r(x;q) [r]_q.
\eqno 2.31
$$
Using 2.25,  the functional equation for $Q_W(x,q)$ may be written as
$$
(1-q/x)  \sum_{r=0}^w A_r(x;q)
\sps x^k(1-xq)\sum_{r=0}^w A_r(1/x;q)
\ses
(1+x^k)(1-q^2)\prod_{i=1}^{k-1}[w_i]_q
\eqno 2.32
$$
While, using 2.27, the  functional equation for  $Q_{W,w}(x;q)$  becomes
\begin{align*}
(1-q/x)\Big( \sum_{r=0}^w & A_r(x;q)\big([r]_qx + q^r[w-r]_q \big) \Big)\sps
\\
\sps x^{k+1}&(1-xq)\Big(\sum_{r=0}^w A_r(1/x;q)\big([r]_q x^{-1}+q^r[w-r]_q\big ) \Big)
\ses
(1+x^{k+1})(1-q^2)[w]_q\prod_{i=1}^{k-1}[w_i]_q
\\
\end{align*}
\vskip -.9in
\hfill $2.33 $

\vskip .5in
\noindent
Notice next that 2.32 
multiplied by $[w]_q$ may be written in the form
\begin{align*}
(1-q/x)\Big( \sum_{r=0}^w &A_r(x;q)
\big([r]_q+q^r[w-r]_q \big) \Big)
\sps
\\
\sps x^k &(1-xq)\Big(\sum_{r=0}^w A_r(1/x;q)\big([r]_q+q^r[w-r]_q\big)  \Big)
=
(1+x^{k})(1-q^2)[w]_q\prod_{i=1}^{k-1}[w_i]_q
\end{align*}

\vskip -.7in
\hfill $2.34 $

\vskip .5in
\noindent
Subtracting 2.34  from 2.33 and dividing by $x-1$  gives
$$
(1-q/x)  \sum_{r=0}^w A_r(x;q)[r]_q
+ x^{k}(1-xq) \sum_{r=0}^w A_r(1/x;q)q^r[w-r]_q \big )  =
x^k (1-q^2)[w]_q\prod_{i=1}^{k-1}[w_i]_q
\eqno 2.35
$$

\newpage
\noindent
This in turn, multiplied by $x$ and  subtracted from 2.33 gives
$$
(1-q/x)  \sum_{r=0}^w A_r(x;q)  q^r[w-r]_q
+ x^k(1-xq) \sum_{r=0}^w A_r(1/x;q) [r]_q
=
(1-q^2)[w]_q\prod_{i=1}^{k-1}[w_i]_q
\eqno 2.36
$$
Using 2.31, 2.36    may also be rewritten as
$$
(1-q/x) A_{W,w}(x;q)
\sps x^k(1-xq)B_{W,w}(1/x;q)
\ses
(1-q^2)[w]_q\prod_{i=1}^{k-1}[w_i]_q
\eqno 2.37
$$
 Likewise  2.35 becomes
$$
(1-q/x) B_{W,w}(x;q)
\sps x^{k}(1-xq)A_{W,w}(1/x;q)   \ses
x^k (1-q^2)[w]_q\prod_{i=1}^{k-1}[w_i]_q
\eqno 2.38
$$
This proves that both 2.22 and 2.23 
are consequences of the functional equations
of $Q_{W}(x;q)$ and $Q_{W,w}(x;q)$.

Conversely suppose that both 2.37 and
2.38 hold true. Then   (using 2.25) we see that their sum is simply
$$
(1-q/x)  Q_W(x;q)[w]_q
\sps x^k(1-xq)Q_W(1/x;q)[w]_q
\ses
(1+x^k) (1-q^2)[w]_q\prod_{i=1}^{k-1}[w_i]_q
$$
and the functional equation of $ Q_W(x;q)$ follows upon division by $[w]_q$.
Similarly, multiplying 2.38 by $x$ and adding it  to 2.37 gives
$$
(1-q/x)  Q_{W,w}(x;q)
\sps x^{k+1}(1-xq)Q_{W,w}(1/x;q)
\ses
(1+x^{k+1}) (1-q^2)[w]_q\prod_{i=1}^{k-1}[w_i]_q
$$
which is the functional equation of $Q_{W,w}(x;q)$.

To complete the proof we  need to show that 2.23 is a consequence of 2.22.
To do this we will use the identities in 2.19 a) and b)   for the schedule $W,w$
in the form
\begin{align*}
&a)\ess\ess\ess A_{W,w}(x;q)\ses q^{d'}x^{k+1}B_{W,w}(1/x;1/q) ,
\cr
&b)\ess\ess\ess  B_{W,w}(1/x;q)\ses q^{d'}x^{-(k+1)}A_{W,w}(x;1/q).
\end{align*}
with  $d'=1+w-1+\sum_{i=1}^{k-1}(w_i-1)$,  the $q$-degree of $P_{W,w}(X_{k+1};q)$. Making these substitutions in 2.22 gives
$$
(1-q/x)  q^{d'}x^{k+1}B_{W,w}(1/x;1/q)
\sps  (1-xq)q^{d'}x^{-1}A_{W,w}(x;1/q)
\ses
(1-q^2)[w]_q\prod_{i=1}^{k-1}[w_i]_q
$$
Replacing $x$ and $q$ by $1/x$ and $1/q$  we get
$$
(1-x/q)  q^{-d'}x^{-k-1)}B_{W,w}(x;q)
\sps  (1-1/xq)q^{-d'}x A_{W,w}(1/x;q)
\ses
-q^{-(d'+1)}(1-q^2)[w]_q\prod_{i=1}^{k-1}[w_i]_q
$$
which  multiplied $q^{d'+1}x^k$ gives
$$
(q/x-1) B_{W,w}(x;q)
\sps x^k (xq-1) A_{W,w}(1/x;q)
\ses
-x^k(1-q^2)[w]_q\prod_{i=1}^{k-1}[w_i]_q
$$
proving 2.23 and completing our proof.
\sas

An immediate Corollary of theorem 4.2 may be stated as follows.
\sas

\noindent{\bol Theorem 2.5}

{\ita If both schedules $W'=(w_1,\ldots,w_{ k-2})$  and $W =(w_1,\ldots,w_{ k-1})$ are tame then the schedule  $W,1=(w_1,\ldots,w_{ k-1},1)$ is also tame as is for any legal schedule
$W"=(W,1,v_2,v_3,\ldots ,v_s)$. In particular, it follows  that from now on there is no loss in assuming that $w_1= 2$  
}

\noindent{\bol Proof}

\newpage
\noindent
We will start by proving  that $W,1$ is tame.
That will be our base case.
 From the hypotheses and Theorem 2.4 it follows that we have the two identities
$$
(1-q/x)  A_{W}(x;q)
\sps x^{k-1}(1-xq)B_{W }(1/x;q)  \ses
  (1-q^2)\prod_{i=1}^{k-1}[w_i]_q
\eqno 2.39
$$
$$
(1-q/x)   B_{W}(x;q)
\sps x^{k-1}(1-xq)A_{W}(1/x;q)  \ses
x^{k-1} (1-q^2)\prod_{i=1}^{k-1}[w_i]_q
\eqno 2.40
$$
Moreover we also have that
$$
P_W(X_k,q)\ses A_W(X_k,q)+x_{k-1}B_W(X_k,q)
$$
since for the schedule $W,1$ the indeterminate $x_k$ only acts on $x_{k-1}$, it follows that
$$
P_{W,1}(X_{k+1};q)\ses
A_W(X_k,q)+x_{k-1}x_kB_W(X_k,q).
\eqno 2.41
$$
In particular we have
$$
Q_{W,1}(x;q)\ses
A_W(x,q)+x^2B_W(x,q).
$$

\sas
\noindent
This given, the functional equation for $Q_{W,1}(x;q)$ must be
\begin{align*}
(1-q/x)\Big(
A_W(x,q)+x^2B_W(x,q)\Big)&\sps
x^{k+1}(1-qx)\Big(
A_W(1/x,q)+x^{-2}B_W(1/x,q)\Big)
=
\\
&=
(1+x^{k+1})(1-q^2)\prod_{i=1}^{k-1}[w_i]_q.
\end{align*}
However, this is easily recognized to be none other than the identity in 2.39 plus the identity in 2.40 multiplied by $x^2$.
\sas

This given, we will show next that we have
$$
P_{W,1,v_2,v_3, \ldots v_s}(X_{k+s},q)= \prod_{i=2}^s[v_i] _q
\Big(A_W(X_{k-1};q)+ x_{k-1}x_{k}x_{k+1}
\cdots x_{k-1+s}
B_W(X_{k-1};q)\Big)
\eqno 2.42
$$
In particular it follows that
\begin{align*}
Q_{W,1,v_2,v_3,\ldots v_s}(x,q)\ses \prod_{i=2}^s[v_i] _q
\Big(A_W(x;q)\sps x^{s+1}B_W(x;q)\Big)
\end{align*}

\vskip -.4in
\hfill
$2.43$ 

\vskip .3in

We start by proving the identity in  2.42 by induction on $s$
with 2.41 as the base case $s=1$.
This given,  let us assume 2.42
true  for $s-1$, that is
$$
P_{W,1,v_2,v_3,\ldots v_{s-1}}(X_{k+s-1},q)
\ses  \prod_{i=2}^{s-1}[v_i] _q
\Big(A_W(X_{k-1};q)\sps x_{k-1}x_{k}x_{k+1}
\cdots x_{k-1+s-1}
B_W(X_{k-1};q)\Big)
$$
Since $x_{k-1+s}$ acts on $x_{k-1+s-1}\cdots x_{k-1+s-v_s} $
and $v_s\le s$ it follows that $x_{k-1} $ is  the last variable which could be acted upon by $x_{k-1+s}$.
This implies that
\begin{align*}
P_{W,1,v_2,v_3,\ldots v_{s}}(X_{k+s},q)
\ses
 \prod_{i=2}^{s-1}[v_i] _q\Big(A_W(X_{k-1};q)[v_s] _q\sps x_{k-1}x_{k}x_{k+1}
\cdots x_{k-1+s-1}x_{k-1+s}[v_s] _q
B_W(X_{k-1};q)\Big)\cr
\end{align*}

\vskip -.2in
\noindent
Proving 2.42 and completing the induction.
Now we may rewrite 2.39  as
$$
(1-q/x)(A_W(x;q) +x^{k+s}(1-qx)(x^{-s-1}B_W(1/x;q))\ses (1-q^2)\prod_{i=1}^{k-1}[w_i]_q
$$
and 2.40 multiplied by $x^{s+1}$ gives
$$
(1-q/x)x^{s+1}  B_W(x;q)+x^{k+s }(1-qx)A_W(1/x;q)\ses x^{k+s } (1-q^2)\prod_{i=1}^{k-1}[w_i]_q
$$
adding these two equalities, multiplying by $ \prod_{i=2}^s[v_i] _q$,
 using 2.43 and setting  $ W'' =W,1,v_2,v_3,\ldots v_s$
 we finally obtain
 $$
 (1-q/x) Q_{  W''  }(x;q) +x^{k+s}(1-qx)(Q_{  W''}(1/x;q))\ses
 (1-q^2)\prod_{i=1}^{k-1}[w_i]_q  \prod_{i=2}^s[v_i] _q(1+x^{k+s} )
$$
proving the functional equation for $Q_{W,1,v_2,v_3,\ldots v_s}(x;q)$ as asserted.
\sa

\centerline {\ita The validity of the Functional equation when all  components of the schedule are equal to  $2$}
\sas

In \cite{thesis} it was discovered that if $W$ is tame
then for certain pairs of integers $u,v$ a legal schedule $uvW$ is also tame.
This given, the following special case will be  essential in our proof of the functional equation,

\sas

\noindent{\bol Theorem  2.6}

{\ita For any legal schedule $W=(w_1,w_2,\ldots ,w_{k-1})$ we have
$$
Q_{2,2,W}(x;q)\ses qx   Q_{W}(x;q)\sps
Q_{2,2,W}(x;q) \Big|_{x}x\sps Q_{2,2,W}(x;q) \Big|_{x^{k+2}}x^{k+2}
\eqno 2.44
$$
In particular, if $W=(w_1,w_2, \cdots ,w_{k-1})$ is
tame then $22W= (2,2,w_1,w_2, \cdots ,w_{k-1})$ is also tame}

\noindent {\bol Proof}

We will establish   2.44 by a bar diagram argument. In the diagrams below, the yellow ``ellipses'' are to represent a  generic bar diagrams for the schedule $W$ with omission of the  initial bars $\UA^{x_{-1}}\DA^{x_0}$ and
$\DA^{x_{-1}}\UA^{x_0}$.
\vskip -.2in
{\begin{figure}[H]
\includegraphics[width=6in]{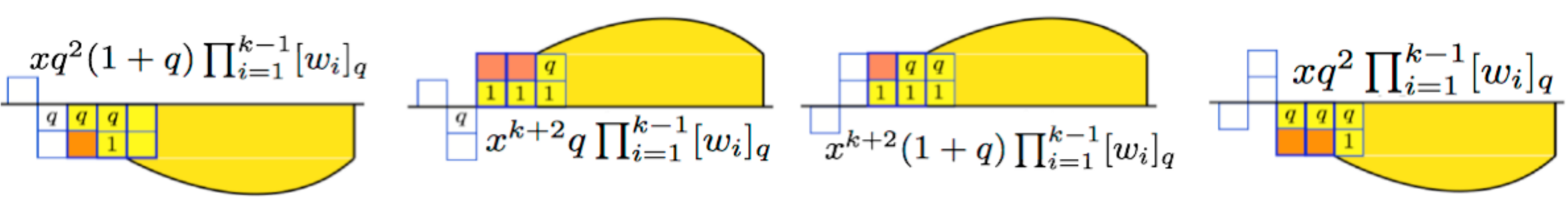}
\end{figure}
\vskip -.28in
{\begin{figure}[H]
\includegraphics[width=4in]{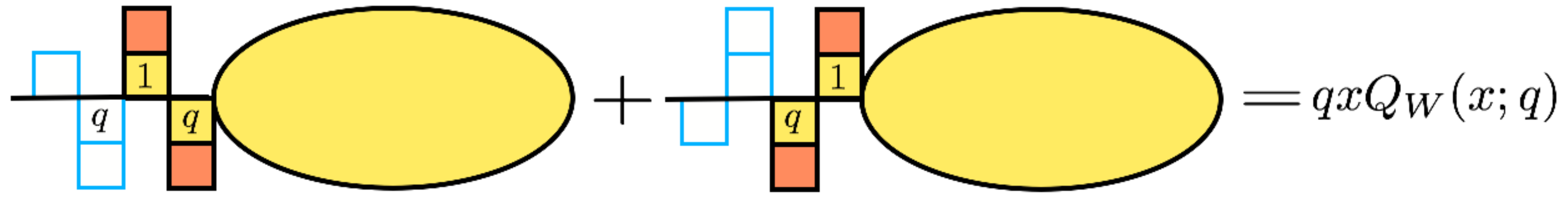}
\end{figure}
\vskip-.16in

The diagrams above will be referred to as $class[1],class[2],\ldots ,class[6]$
listed  as we scan  them from left to right and from top to bottom.
We will use  them to represent a decomposition of the collection of the bar diagrams
of the schedule $2,2,W$ into ``classes '' according to the positions of
the first four bars. That is  can simply represent them by the symbols
\begin{align*}
\UA^{x_{-1}}\DA^{x_0}\DA^{x_1}\DA^{x_2},\ess\ess
&\UA^{x_{-1}}\DA^{x_0}\UA^{x_1}\UA^{x_2},\ess\ess
\ess\ess\ess\ess\ess
\DA^{x_{-1}}\UA^{x_0}\UA^{x_1}\UA^{x_2}\ess\ess
\DA^{x_{-1}}\UA^{x_0}\DA^{x_1}\DA^{x_2}
\\
&\UA^{x_{-1}}\DA^{x_0}\UA^{x_{1}}\DA^{x_2},
\ess\ess\ess\ess\ess\ess\ess
\DA^{x_{-1}}\UA^{x_0}\DA^{x_{1}}\UA^{x_2}.
\end{align*}
In $class[1]$ and $ class[3]$  bar $x_2$ has no choices, being only of length $2$ in each case it has to be on the same side of the ground line as bar $x_1$. For $class[2]$ and $class[4]$ we chose to depict bar $x_2$ on the same side of bar $x_1$.
This forces all the remaining bars to be  totally yellow. We can clearly see that for $class[1]$,  $class[2]$, $class[3]$  and $class[4]$,
the   polynomial above or below each diagram
gives precisely the contribution that the class makes
to $Q_{2,2,W}(x;q)$.

Finally, we see that $ class[5]$
and $class[6]$ are obtained by choosing the other alternative for bar $x_2$.
The second author's  beautiful idea in her thesis, is to view $class[5]$ and $class[6]$,
as the result of prepending the pair
$
\UA^{x_{-1}}\DA^{x_0}
$
to the collection of bars diagrams of the schedule $W$ that start with
$
\UA^{x_{-1}}\DA^{x_0}
$
and respectively prepending the pair
$
\DA^{x_{-1}}\UA^{x_0}
$
to the collection of bars diagrams of the schedule $W$ that start with
$
\DA^{x_{-1}}\UA^{x_0}
$. A look at the initial required labelings reveals that the contribution of these two classes to the polynomial
$Q_{2,2,W}(x;q)$ is none other than $xq Q_W(x;q)$. This given, since we can easily see that $class[1]$ and $class[4]$
yield the coefficient of $x$ and
 $class[2]$ and $class[3]$ yield the coefficient of $x^{k+1}$ we obtain  the equalities
 \vskip -.2in
$$
a)\ess
Q_{2,2,W}(x;q) \Big|_{x}=
(2q^2+q^3)\prod_{i=1}^{k-1}[w_i]_q
,\ess\ess\ess\ess\ess\ess
\ess
b)\ess\ess
Q_{2,2,W}(x;q) \Big|_{x^{k+2}}=
(1+2q)\prod_{i=1}^{k-1}[w_i]_q
\eqno 2.45
$$
This completes our proof of 2.44.

\newpage

We are left with checking that the polynomial $Q_{2,2,W}(x;q)$ satisfies the
functional equation. For convenience let us express it as
$\ess
Q_{2,2,W}(x;q)=b_1(q)x+xqQ_w(x;q)+b_{k+2}(q)x^{k+2}
$.
This gives
\begin{align*}
(1-q/x)Q_{2,2,W}&(x;q)+x^{k+2}(1-qx)Q_{2,2,W}(1/x;q)\ses
\\
&=b_1x-qb_1
+b_{k+2}x^{k+2}-qb_{k+2}x^{ k+1}
+ b_1x^{ k+1}-qb_1x^{k+2}+ b_{k+2}-qb_{k+2}x\cr
&
\bigsp\sps
xq(1-q/x)Q_W(x;q)+x^{k+2}(1-qx)(q/x)Q_W(1/x;q)
\\
&= (b_{k+2} -qb_1)(1+x^{k+2})
\sps (b_1 -qb_{k+2})(x
+x^{ k+1}) \sps
\cr
&\bigsp\bigsp
\hbox{(since $W$ is tame)}
\ess\ess\sps qx(1+x^k)(1-q^2) \tttt{\prod_{i=1}^{k-1}[w_i]_q}
\\
\end{align*}

\vskip -.18 in
\noindent
But from 2.45 it follows that
$$
(b_1 -qb_{k+2})\sps q(1-q^2) {\prod_{i=1}^{k-1}[w_i]_q}\ses
(q^3-q){\prod_{i=1}^{k-1}[w_i]_q}
\sps q(1-q^2){\prod_{i=1}^{k-1}[w_i]_q}\ses
0.
$$
Finally, we are left with
$$
b_{k+2} -qb_1=
\Big((1+2q) -q(2q^2+q^3)\Big)
\prod_{i=1}^{k-1}[w_i]_q
=
(1-q^2)[2]_q[2]_q
\prod_{i=1}^{k-1}[w_i]_q.
$$
This completes our argument.
\sa

Since the schedules $W=\{\emptyset\}$  and  $W=(2)$ are clearly tame
we have the following important corollary.
\sas

\noindent{\bol Theorem 2.7}

{\ita The schedules of the form $W=(2,2,2,\ldots ,2)$  are   tame.}
\sas

We should mention that the first proof of this result appeared in \cite{DAGA}, however the present proof proves much more with considerably less effort.
\sa

\section{Proof of the Functional equation}

After  so many years of wondering about the validity of this result, the simplicity and modality of its proof is stunning. The proof assumes that all legal schedules
$W=(w_1,w_2,\ldots,w_{h-1})$ with $h<k$ are tame. Then proves it for $k$ by recursing over a finite set.

More precisely, let $\CF[k]$ denote the family of all legal schedules $W=(2,w_2,\ldots w_{k-1})$ whose components are all $\ge 2$, totally ordered by  lex order.
The minimal element of $\CF[k]$ is therefore the schedule   $W=(2, 2,\ldots 2)$ with $k-1$ $2's$. The proof constructs
two recursions $\phi(W)=[W';W'']$ and
$\psi(W)=[W';W'';W''']$ with the following properties
\begin{enumerate}
\item   $W',W''<_{lex}W$, for  $\phi$ and $W',W'',W'''<_{lex}W$  for $\psi$.
\item   For each $W$ in $\CF[k]$ one and only one of $\phi$ or $\psi$ applies.
\item   For $\phi$: $Q_W=(1+q)Q_{W'}-q Q_{W''}$
\item      For $\psi$: $Q_W=\tttt{q\over 1+q}Q_{W'}+Q_{W''}-\tttt{q\over 1+q} Q_{W'''}$
\item   The only ones of
$W',W''<_{lex}W$, for  $\phi$ and $W',W'',W'''<_{lex}W$  for $\psi$,
that may be not in $\CF[k]$ will have one of their component equal to $1$.
\end{enumerate}

\noindent
It is easily shown that these recursions
reduce the tameness of $W$ to the tameness of their predecessors $W',W'' $ or $W',W'',W'''$ as the case may be. The base cases turn out to be legal schedules
$W=(w_1,w_2,\ldots ,w_{k-1})$ with at least one component $1$ or the schedule
$W=(2, 2,\ldots 2)$ with $k-1$ $2's$.
Since the base cases have been shown to be tame it follows that all elements of  $\CF[k]$ are tame.

The remainder of this section is dedicated
to the construction of the recursions $\phi$ and $\psi$ and proving their stated properties. We will start with helpful   definitions and some auxiliary observations.

Recall that $W=(w_1,w_2,\ldots ,w_{k-1})$ is legal
if $w_{i }\le w_{i-1}+1$ for all $2\le  i\le k-1$.
It is clear that all schedules $W\in \CF[k]$, except the first,
have at least one increase. The last increase will be called the {\ita canonical increase}  of $W$. If $w_{h-1}=v$ and $w_h=v+1$ is an increase of $W$ we will say that bar $x_{r}$ {\ita splits} $w_{h-1}$ and $w_h$ if
$r-w_r=h$. Recalling  that bar $x_r$ acts on bars $x_{r-1},x_{r-2},\ldots,x_{r-w_r}$. That means that bar $x_r$ acts on bar $x_h$ but not on bar $x_{h-1}$. If no bar splits the canonical increase of $W$ we will write $W=U,v,v+1,V$ where $v,v+1$ is the canonical increasing pair and define
\vskip-.13in
$$
\phi(W)=[U,v,v,V\,\, ;\, \,U,v-1,v,V]=[W'\, ;\, W'']
\eqno 3.1
$$
If some bar does split the canonical increase then let $x_r=x_{h+j}$ be the first bar that does. Since that means we
have $h+j-w_{h+j}=h$ we can
write $W=U,v,v+1,V',j,V''$ with $j$  in position $h+j$. This given, we define
\begin{align*}
\psi(W)=[
U,v,v+1,V',j-1,V''\,\, &;\,\,
U,v,v,V',j+1,V''\,\, ;\\
&;
U,v-1,v ,V',j+1,V'']=[W'\, ;\, W''\, ;\, W''']
\end{align*}

\vskip -.35in
\hfill
$
3.2
$
\vskip .2in
\noindent
It is clear from 3.1 and 3.2 that property
(1) of $\phi$ and $\psi$ is satisfied.
The following two propositions show that both are well defined.
\sas

\noindent{\bol Proposition 3.1}

{\ita The schedules $W'$ and $W''$ in 3.1 are legal}

\noindent{\bol Proof}

The only problem arises if $V$ would start with a component $v+2$. But if that were the case then we could write
$W=U,v,v+1,v+2,V'$ and contradict that
$v,v+1$ is the last  increase of $W$.
\sas

\noindent{\bol Proposition 3.2}

{\ita The schedules $W',W''$ and $W'''$ in 3.2 are legal}

\noindent{\bol Proof}

Notice first that $W''$ and $W'''$ would not be legal if $V'$ could start with
$v+2$. But then again as we have seen
in the previous proof, that would contradict that $v,v+1$ is the last increase of $W$. Next,  the legality
of $W'$  would not hold if $V''$
could start with $j+1$. But if that were the case then $w_{h+j }=j  $ and $w_{h+j+1}=j+1$ would again contradict that the  pair $v,v+1$ is canonical for $W$. The last remaining issue is the $j+1$ in $W''$
if $V'$ could end with a component less than $j$. However the
legality of $W$ would force that component to be $j-1$ and we could write
$W=U,v,v+1,V',j-1,j,V''$ again  contradicting
that the pair $v,v+1$ is canonical for $W$.
\sas
\noindent{\bol Remark 3.1}

For future purposes we need to focus on the schedule obtained from
$W=U,v,v+1,V',j,V"$ by setting
$\widetilde W=U,v,v+1,V',j+1,V"$.
Notice first that since $W$  is legal  the last component of $V'$ must be $\ge j-1$. But equality here cannot hold for two reasons. Firstly, since the pair  $j-1,j$ would contradict that $v,v+1$
is canonical for $W$. Secondly, since we must recall that   $x_r=x_{h+j}$ was chosen to be the first bar that splits the pair $v,v+1$. This shows that
$\widetilde W$ is legal. Nevertheless the last component in $V'$ could very well be $j$ and the pair $v,v+1$ would cease to be canonical for $\widetilde W$. However, what is of crucial importance is that there  are no bars in  $\widetilde W$ that split  the pair $v,v+1$.  The reason is that $j$ in $W$ was chosen to be the length of first bar that did. Now in $\widetilde W$ that length is changed to $j+1$. Could there be in $ V''$ a bar
$x_{r'}$ that splits $v,v+1$? The answer in no. Let us recall that  the position of $v+1$ in $W$ is $h$, thus for $x_{r'}$ to act on $v+1$ and not on $v-1$ we must have
$r'-w_{r'}=h $. This gives $w_{r'}=r'-h$.
Since $r'>r=h+j$ then  we would also have
$
w_{r'}> h+j-h=j=w_r
$
But that cannot happen in $W$ since the last increase is $v,v+1$.
\sas

Our definitions of $\phi$ and $\psi$ guarantee  property $(2)$. Property (5)
is assured by 3.1 and 3.2. In fact, if $W$ has  any component $1$ there is no need to apply $\phi$ or $\psi$ to it.
This given, we see from 3.1 that an application of $\phi$ may produce a $W''$ with a component equal to $1$ by the replacement $v\RA v-1$. Similarly we see from 3.2 that $W'$ by  $j\RA j-1$ and $W'''$ by $v\RA v-1$ may end up with  a component $1$.

Thus we are left with proving $(3)$ and $(4)$ and deriving from these identities that the tameness of the corresponding  $W \, 's$ forces the tameness of $W$.
The following result not only proves $(3)$, but also plays a role in the  proof of $(5)$.

\newpage
\noindent{\bol Proposition 3.3}

{\ita Let $W=U,v,v+1,V$, where $v,v+1$ is not necessarily the last increase,  but suppose that there is  no bar $x_r$ that splits the pair $v,v+1$ then we have}
$$
Q_{U,v,v+1,V}[x,q]
=(1+q)Q_{U,v,v,V}[x,q]-qQ_{U,v-1,v,V}[x,q].
\eqno 3.3
$$
\noindent{\bol Proof}

Let $v,v+1$ be the components $w_{h-1},w_{h}$
of $W$. As before let us denote by $W'$ and $W''$ the schedules appearing in the right hand side of 3.3. If $BD,BD',BD''$ are corresponding (same up-down)  bar diagrams of $ W, W' , W'' $ respectively, it will be convenient to use the short  notation
$$
 \Gamma^{W-[2]W'+qW''} (BD) = \Gamma^{W}(BD)-[2]\Gamma^{W'}(BD')+q\Gamma^{W''}(BD''),
 \eqno 3.4
 $$
where $\GG^W,\GG^{W'},\GG^{W''}$ denote taking the weights contributed by these bar diagrams to the respective polynomials
$P_W,P_{W'},P_{W''}$. Now let $BD^U$ and $BD^V$ be the portions of $BD$ contributed by $U$ and $V$ respectively. Our proof of 3.3
is based on establishing the following identities.
\begin{align*}
 \Gamma^{W-[2]W'+qW''}(BD^U \downarrow \downarrow  BD^V )&=\Gamma^{W-[2]W'+qW''}(BD^U \uparrow \uparrow  BD^V)=0,
 \cr
  \Gamma^{W-[2]W'+qW''}(BD^U \downarrow \uparrow  BD^V)&=-\Gamma^{W-[2]W'+qW''}(BD^U \uparrow \downarrow  BD^V).
\end{align*}

\vskip -.33in
\hfill
$
 3.5
$
\vskip .15in

\noindent
where the symbol ``$BD^U \downarrow \downarrow  BD^V$'' denotes the bar diagram $BD$ which starts with $BD^U$ followed with bars $x_{h-1},x_h$ both down and finishing with $BD^V$.
In  the same vein as  in 3.4, we are also  requiring here the  Bar Diagrams, $BD'$ and $BD''$ to have the same up or down state of bars $x_{h-1},x_h$ as in  $BD$.
To prove these identities let us recall that
the contribution to the weight of   bar $x_j$ in a diagram $BD$ of $W$  is $[a_j]_q$ if  bar $x_j$ is  up  and
$[w_j]_q -[a_j]_q=q^{a_j}[w_j-a_j]_q$
if bar $x_j$ is  down.  Where we have denoted by $a_j$ the number of up-arrows acted upon by bar $x_j$ when it is up.

Since the three schedules $W$,$W'$, $W''$
only differ at the two entries in positions $h-1$ and $h$, every bar diagram $BD$, $BD'$, $BD''$  shares the same $a_j$ sequence except for  $a_{h-1}$ and $a_h$. Keeping this in mind we only need to distinguish two cases:
\sas

{\bf Case 1:} $BD$ has bar $x_{h-v-1}$ down. There are four states of bars $x_{h-1}$ and $x_h$ as listed in the following table. (Here, next to each arrow we place the number of up bars it acts upon, when that bar is up)
$$
\begin{matrix}
&|&w_{h-1}= v &  w_h=v+1 &|&  w'_{h-1}=v & w'_{h}=v &|& w_{h-1}''=v-1 & w_h''=v
 \\
     x_{h-v-1}\DA &|& \DA b & \DA b &|& \DA b & \DA b &|& \DA b & \DA b
\\
     x_{h-v-1}\DA &|& \UA b & \UA b+1 &|& \UA b & \UA b+1 &|& \UA b & \UA b+1
\\
     x_{h-v-1}\DA &|& \DA b & \UA b &| & \DA b &  \UA b &|& \DA b & \UA b
\\
     x_{h-v-1}\DA &|& \UA b & \DA b+1 &|& \UA b &  \DA b+1 &|& \UA b & \DA b+1
\\
\end{matrix}
$$
Thus we have
\begin{align*}
\Gamma  ^{W-[2]W'+qW''}& (BD^U\da\da BD^V)=\cdots \times\cr
& \times\Big(([v]-[b])([v+1]-[b])-(1+q)([v]-[b])^2+q([v-1]-[b])([v]-[b])\Big)=0,\\
\Gamma  ^{W-[2]W'+qW''}&(BD^U\ua\ua BD^V) =\cdots \cr
&\times\Big([b][b+1]-(1+q)[b][b+1]+q[b][b+1]\Big)=0,\\
\Gamma ^{W-[2]W'+qW''}&(BD^U\da\ua BD^V)=\cdots\times \\\
& \times \Big(([v]-[b])[b]-(1+q) ([v]-[b])[b] +q ([v-1]-[b])[b]\Big)=-q^v[b]\cdots,\cr
\Gamma ^{W-[2]W'+qW''}&(BD^U\da\ua BD^V)=\cdots\times  \\
&\times \Big([b]([v+1]-[b+1])-(1+q)[b]([v]-[b+1])+q [b] ([v]-[b+1])\Big)=q^v[b]\cdots.
\\
\end{align*}

\sas
{\bf Case 2:} $BD$ has $x_{h-v-1}\uparrow$. There are also four states of bars $h-1$ and $h$ as listed in the following table.
$$
\begin{matrix}
      &||& w_{h-1}= v & w_h=v+1 &||& w'_{h-1}=v & w'_{h}=v &||& w_{h-1}''=v-1 & w_h''=v \\
     x_{h-v-1}\uparrow &||&  \downarrow b+1 & \downarrow b+1 &||&  \downarrow b+1 & \downarrow b  &||&  \downarrow b & \downarrow b \\
     x_{h-v-1}\uparrow &||&  \uparrow b+1 & \uparrow b+2 &||& \uparrow b+1 & \uparrow b+1 &||& \uparrow b &   \uparrow b+1 \\
     x_{h-v-1}\uparrow &||&  \downarrow b+1 & \uparrow b+1 &||& \downarrow b+1 & \uparrow b &||& \downarrow b &   \uparrow b \\
     x_{h-v-1}\uparrow &||&  \uparrow b+1 & \downarrow b+2 &||& \uparrow b+1 & \downarrow b+1 &||& \uparrow b & \downarrow b+1
\end{matrix}
$$
Thus we have
\begin{align*}
\Gamma& ^{W-[2]W'+qW''}(BD^U\da\da BD^V)=\cdots\times \\\
&\times\Big(([v]-[b+1])([v+1]-[b+1])-(1+q)([v]-[b+1])([v]-[b])+q([v-1]-[b])([v]-[b])\Big)=0,\\
\Gamma& ^{W-[2]W'+qW''}(BD^U\ua\ua BD^V) =\cdots \\
&\times \Big([b+1][b+2]-(1+q)[b+1][b+1]+q[b][b+1]\Big )=0,\\
\Gamma&^{W-[2]W'+qW''}(BD^U\da\ua BD^V) =\cdots\times\\
& \times\Big(([v]-[b+1])[b+1]-(1+q) ([v]-[b+1])[b] +q ([v-1]-[b])[b] \Big) =([v]-[b+1])\cdots,\\
\Gamma&^{W-[2]W'+qW''}(BD^U\da\ua BD^V)=\cdots\times \\
&=\times \Big([b+1]([v+1]-[b+2])-(1+q)[b+1]([v]-[b+1])+q [b] ([v]-[b+1])\Big)=-([v]-[b+1])\cdots.
\end{align*}
Summing the weights of bar diagrams considered in these two cases proves 3.5 and completes the proof of the proposition.
\sas

Let $W=U,v,v+1,V$ where $v,v+1$ is the canonical increasing pair of $W$ and suppose that there is  a bar $x_r$ that splits the pair $v,v+1$. As we did before, let $r$ be the smallest one. That is $r=h+j$ with $j$ the smallest such that $h+j-w_{h+j}=h$. Thus we can  write
$W=U,v,v+1,V',j, V''$ with
$j$ in position $h+j$.
\sas

 This given, from Remark 3.1 it follows that the schedule
$\widetilde W=U,v,v+1,V',j+1, V''$ is not only legal but Proposition 3.3 is
also applicable to it. This gives the identity
$$
  Q_{\widetilde{W}}[x,q] -(1+q)  Q_{U,v,v,V',j+1,V''}[x,q]+  qQ_{U,v-1,v,V',j+1,V''}=0
\eqno 3.8
$$
\sas

We claim that we can also prove
\sas

\noindent{\bol Proposition 3.4}

{\ita Let the pair  
$v,v+1$ be canonical for the schedule $W=U,v,v+1,V',j,V''$ with $v+1$ in position $h$ and $j $ in position
$r=h+j$, so that  bar $x_r=x_{h+j}$
splits $v,v+1$, we can still prove
the identity
$$
 Q_{U,v,v+1,V',j+1,V''}[x,q] -(1+q)  Q_{U,v,v+1,V',j,V''}[x,q]+  qQ_{U,v,v+1,V',j-1,V''}[x,q]=0
\eqno 3.9
$$
In particular  it follows that we have}
$$
Q_{W} =  \tttt{q\over1+q}Q_{U,v,v+1,V',j-1,V''}
 + Q_{U,v,v,V',j+1,V''} -  \tttt{q\over 1+q}Q_{U,v-1,v,V',j+1,V''}.
\eqno 3.10
$$
\noindent{\bol Proof}
\def \TW  {{\widetilde W}}
\def \TV  {{\widetilde V}}
\def \TW  {{\widetilde W}}
For convenience let us write
$$
\TW=U,v,v+1,\TV\ess\ess\ess
W'=U,v,v+1,\TV'\ess\ess\ess
W''=U,v,v+1,\TV''\ess\ess\ess
$$
with
$$
\WV= V',j+1,V'',    \ess\ess\ess\ess
\WV'=V',j ,V'',    \ess\ess\ess\ess
\WV''= V',j-1,V'',
$$
This given, we will  prove 3.9 by simultaneously summing the weights
$$
  \GG^\TW \big( BD^{ U}\DA\DA BD^\WV\big),
\ess\ess\ess\ess\ess\ess\ess\ess
  \GG^{W'} \big( BD^{ U}\DA\DA BD^{\WV'}\big),
\ess\ess\ess\ess\ess\ess\ess\ess
 \GG^{W''} \big( BD^{ U}\DA\DA BD^{\WV''}\big),
\eqno 3.11
$$
then doing the same by the two arrows $``\DA\DA$''
respectively replaced by ``$\UA\UA$'', ``$\DA\UA$'' and ``$\UA\DA$''.

\newpage
Notice that the uniqueness of a splitter of the  canonical increasing  pair assures
(see the end of Remark 3.1) that  there will be no bar $x_{h+i}$ with  $i\neq j$ in $\WV=V',j+1, V''$, $\WV=V',j,V''$ and $\WV=V',j-1,V''$
that will also split $v,v+1$.
Consequently, we  only  need to focus on the $a_i$-sequences in positions $h-1,h$ and $h+j$.

Keeping this in mind,  we will verify the following equalities. The identity in 3.9  then follows by summing over all $BD^U$ and $B^\WV$ with  the four  up,down states of bars $x_{h-1}$ and $x_h$.
\begin{align*}
  a)\ess\ess \Gamma^{\TW-[2]W'+qW''}(BD^U\downarrow \da BD^\WV)&=\Gamma^{\TW-[2]W'+qW''}(BD^U\uparrow \uparrow  BD^\WV)=0,\\
  b)\ess\ess  \Gamma^{\TW-[2]W'+qW''}(BD^U\uparrow \downarrow  BD^\WV)&=-\Gamma^{\TW-[2]W'+qW''}(BD^U\downarrow \uparrow  BD^\WV).
\end{align*}

\vskip -.35in
\hfill
$
3.12
$

\vskip .2in

\noindent
To prove these three identities, we observe that $\TW,W',W''$ only differ at the $h+j$-th entry, so does their corresponding weight of any particular bar diagram. We also need to consider the state of bars $x_{h-1}$ and $x_h$ since they cause the difference of the weights of bar $h+j$ with respect to $\TW,W',W''$.

We use the following table for the cases when bar $h+j$ is up or down:
$$
\begin{matrix}
    w_{h-1}=v & w_{h}=v+1 &||& w_{h+j}=j+1 & w_{h+j}=j & w_{h+j}=j-1 \\
    \downarrow c & \downarrow c&||& \updownarrow a & \updownarrow a & \updownarrow a \\
    \uparrow c & \uparrow c+1&||& \updownarrow a+2 & \updownarrow a+1 & \updownarrow a \\
    \downarrow c & \uparrow c& ||&\updownarrow a+1 & \updownarrow a+1 & \updownarrow a \\
    \uparrow c & \downarrow c+1&||& \updownarrow a+1 & \updownarrow a & \updownarrow a
\\
\end{matrix}
$$
The case when bar $x_{h+j}$ is down can be derived  from the case when bar $x_{h+j}$ is up. In fact, since $[j+1]-(1+q)[j]+q[j-1]=0$, the contribution of bar $x_{h+j}$ down is simply that of bar $x_{h+j}$ up times $-1$. Thus we assume that bar $x_{h+j}$ is up in the computations below.
\sas

{\bf Case 1:} Bars $x_{h-1}$ and $x_{h}$ are both down or both up. We only need to concentrate on bar $x_{h+j}$ since its contribution forces the whole sum to vanish. In fact, have
\begin{align*}
\Gamma^{\TW-[2]W'+qW''}(BD^U\da\da BD^\WV) &=\cdots
([a]x-(1+q)[a]x+q[a]x) =0,\\
  \Gamma^{\TW-[2]W'+qW''}(BD^U\ua\ua BD^\WV) &=\cdots
([a+2]x-(1+q)[a+1]x+q[a]x) =0.
\\
\end{align*}
\vskip -.2in

{\bf Case 2:} Exactly one of bars $x_{h-1}$ and $x_h$ is up. In this  case  we need to concentrate on all three bars $x_{h-1}$, $x_h$ and $x_{h+j}$. We have

\begin{align*}
\Gamma^{\TW-[2]W'+qW''}(BD^U\da\ua BD^\WV) &=\cdots
q^c[v-c][c]\Big([a+1]x-(1+q)[a+1]x+q[a]x\Big)\\
&=-q^{a+c+1}[c][v-c]x\cdots,\\
  \Gamma^{\TW-[2]W'+qW''}(BD^U\ua\da BD^\WV) &=\cdots
[c]q^{c+1}[v+1-c-1]\Big([a+1]x-(1+q)[a]x+q[a]x\Big)\\\
&=q^{a+c+1}[c][v-c]x \cdots.
\\
\end{align*}
This completes the proof of the identity in 3.9.
\sas

Setting the  two identities 3.8 and 3.9 side by side we notice that we can write them in the form
\vskip -.18in
\begin{align*}
  Q_{\widetilde{W}}  &\ses (1+q)  Q_{U,v,v,V',j+1,V''} -  qQ_{U,v-1,v,V',j+1,V''}
\\
  Q_{\widetilde{W}}  &\ses (1+q)  Q_W -  qQ_{U,v,v+1,V',j-1,V''}
\\
\end{align*}
\vskip -.18in
\noindent
Eliminating  $ Q_{\widetilde{W}}$ proves the identity
$$
Q_{W} =  \tttt{q\over1+q}Q_{U,v,v+1,V',j-1,V''}
 + Q_{U,v,v,V',j+1,V''} -  \tttt{q\over 1+q}Q_{U,v-1,v,V',j+1,V''}.
$$
which we easily recognise as 3.10. 
This completes our proof of Proposition 3.4.
\sas

The next and final step is to prove the following beautiful  fact

\newpage
\noindent{\bol Proposition 3.5}

{\ita Let $W=(w_1,w_2,\ldots,w_{k-1})$ be a legal schedule
\begin{align*}
& \hbox{$\bf (a)$ If we use $\phi(W)=(W',W'')$
and   $W'$, $W''$ are tame  then $W$ is tame,}
\\
&
\hbox{$\bf (b)$  If  we use  $\psi(W)=(W',W'',W''')$ and $W'$, $W''$, $W'''$ are tame then $W$ is tame.}
\end{align*}}
\noindent{\bol Proof}

Recall that a legal schedule $W=(w_1,w_2,\ldots ,w_{k-1})$ is tame if an only if
$$
(1- {q/ x})Q_W(x;q)\sps x^k(1-qx)Q_W(1/x;q)
\ses (1+x^k)(1-q^2)\prod_{i=1}^{k-1}[w_i]_q
\eqno 3.13
$$
This given, notice first that the $\BQ[q]$-linear operator
$$
\CL F(x;q)=(1- {q/ x})F(x;q)\sps x^k(1-qx)F(1/x;q)
$$
sends a Laurent polynomial in $x$ 
into another such polynomial. For convenience, for any positive integral vector $U=(u_1,u_2,\ldots,u_r)$  set
$
\Pi[U]=\prod_{i=1}[u_i]_q
$. With this notation we may write 3.13 as
$$
\CL Q_W\ses (1+x^k)(1-q^2)\Pi[W]
$$
\sas

If  we use $\phi$ then from Proposition 3.3 
and $W=U,v,v+1,V$  then it follows that
 $$
Q_{U,v,v+1,V}[x,q]
=(1+q)Q_{U,v,v,V}[x,q]-qQ_{U,v-1,v,V}[x,q].
\eqno 3.14
$$
Thus   the $\BQ[q]$-linearity of $\CL$ gives
$$
\CL Q_{U,v,v+1,V}= (1+q)\CL Q_{U,v,v,V}-q\CL Q_{U,v-1,v,V}
\eqno 3.15
$$
If $W'=U,v,v,V$ and $W''=U,v-1,v,V$ are tame then
3.15 becomes
$$
\CL Q_{W}= (1+x^k)(1-q^2)
\Pi[U]\Big((1+q)[v]_q[v]_q-q[v-1]_q[v]_q
\Big)\Pi[V]
$$
and the tameness of $W$ is reduced to showing the  trivial identity
$$
[v]_q[v+1]_q\ses (1+q)[v]_q[v]_q-q[v-1]_q[v]_q.
$$
\sas

If  we use $\psi$ and $W=U,v,v+1,V',j, V''$
then  from Proposition 3.4 it follows that
$$
Q_{W} =  \tttt{q\over 1+q}Q_{U,v,v+1,V',j-1,V''}
 + Q_{U,v,v,V',j+1,V''} -  \tttt{q\over 1+q}Q_{U,v-1,v,V',j+1,V''}.
\eqno 3.16
$$
Thus the tameness of $W'$,$W''$ and  $W'''$ gives
$$
\CL Q_{W}= (1+x^k)(1-q^2)
\Pi[U][v]_q
\Big(
\tttt{q\over 1+q}[v+1]_q[j-1]_q+
[v]_q[j+1]_q-\tttt{q\over 1+q}[v-1]_q[j+1]_q
\Big)\Pi[V']\Pi[V'']
$$
and the tameness of $W$ is reduced to showing the identity
\begin{align*}
[v+1]_q[j]_q
&\ses \tttt{q\over 1+q}[v+1]_q[j-1]_q+
[v]_q[j+1]_q-\tttt{q\over 1+q}[v-1]_q[j+1]_q
\\
&\ses [j]_q\sms [v]_q+[v]_q[j+1]_q
\ses  [j]_q  +[v]_q q [j]_q
\end{align*}
which is trivially true.
This completes the proof of Proposition 3.5 and the proof
of the tameness of all legal schedules.
\sas

In the display below, our recursion is applied to the 14 schedules of length $4$ in lex order, (omitting  $[2,2,2,2]$).  As we can plainly see, by reading the columns from top to bottom and from left to right,
every schedule is followed by its image by $\phi$ or $\psi$ as dictated by the algorithm.

\vskip-.15in
{\begin{figure}[H]
\includegraphics[width=5.5in]{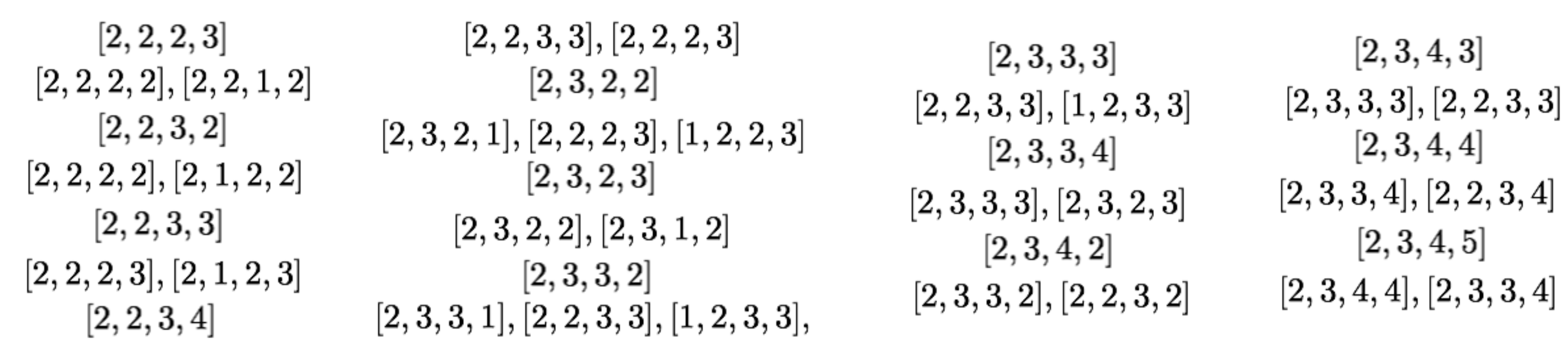}
\end{figure}
 \vskip -.38in

\section{ From Bar Diagrams to Parking Functions: open problems and conjectures }

We will start by making clearer the connection between Parking Functions and Labeled Bar Diagrams. This may require the repetition of some material covered in the first section. However it is very important that this connection is clearly understood for
future applications of the results of this paper.

To carry this out we need first a detailed description of an algorithm used in \cite{thesis} to construct all the Parking Functions with prescribed diagonal cars. This algorithm extends to the case of Frobenius series what Haglund and Loehr did in \cite{HAGLO} for Hilbert series of Diagonal Harmonics.

 This is better understood by means of a specific example. Given the permutation
$$
\sig\ses [5, 6, 2, 4, 7, 8, 1, 3]
\eqno 4.1
$$
The first step is to break it into increasing runs
$$
runs(\sig)\ses 5, 6\, |\,  2, 4, 7, 8\, |\,1, 3
\eqno 4.2
$$
Our task is to construct all the Parking Functions with sets of cars $\{1,3\}$, $\{2, 4, 7, 8\}$, $\{ 5, 6\}$, in diagonals $0,1,2$  respectively.
The display below exhibits such a Parking Function
\vskip -.1in
{\begin{figure}[H]
\includegraphics[width=3.6in]{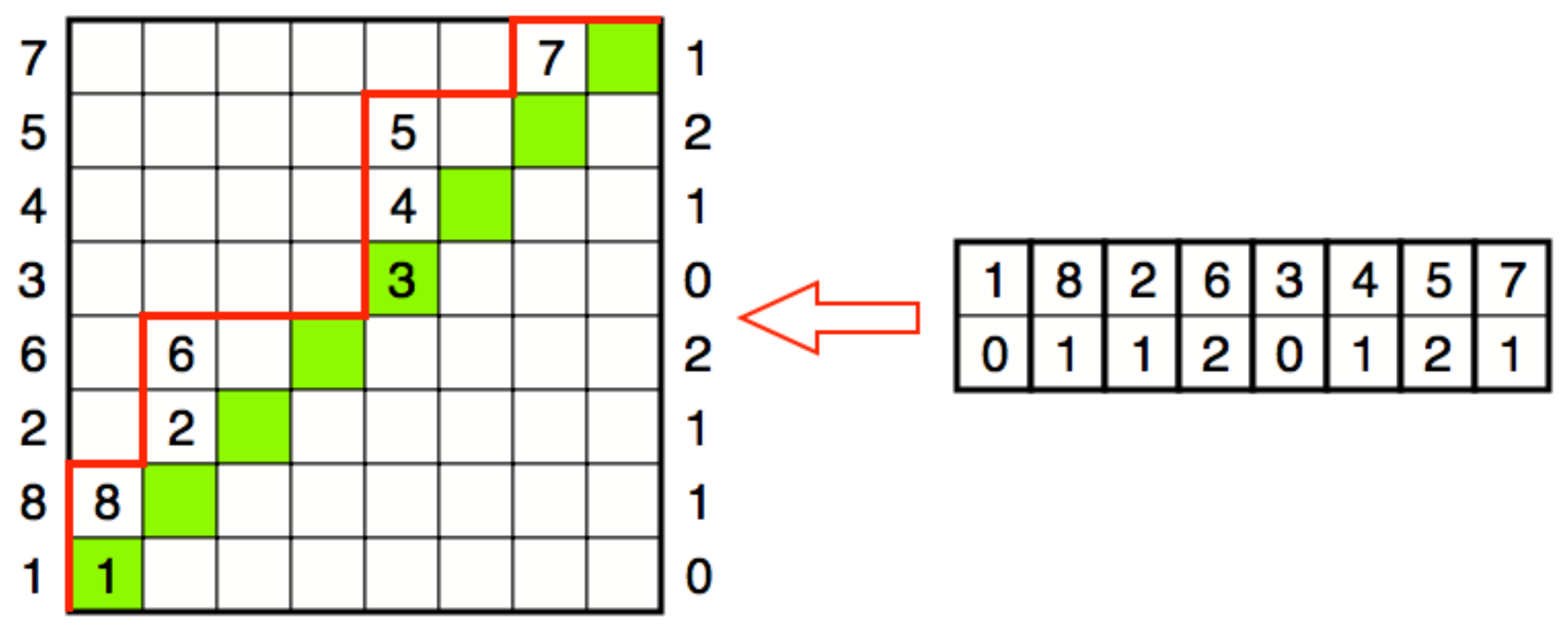}
\end{figure}
\vskip -1in
\hfill
$4.3$
\vskip  .7in
\noindent
 The idea is to consider a car $c$ to be placed in diagonal $d$ as a vertical domino labeled by  $c$ over $d$.
In particular the runs of $\sig$ are thus converted into the three runs of dominoes displayed below
\vskip -.1in
{\begin{figure}[H]
\includegraphics[width=2in]{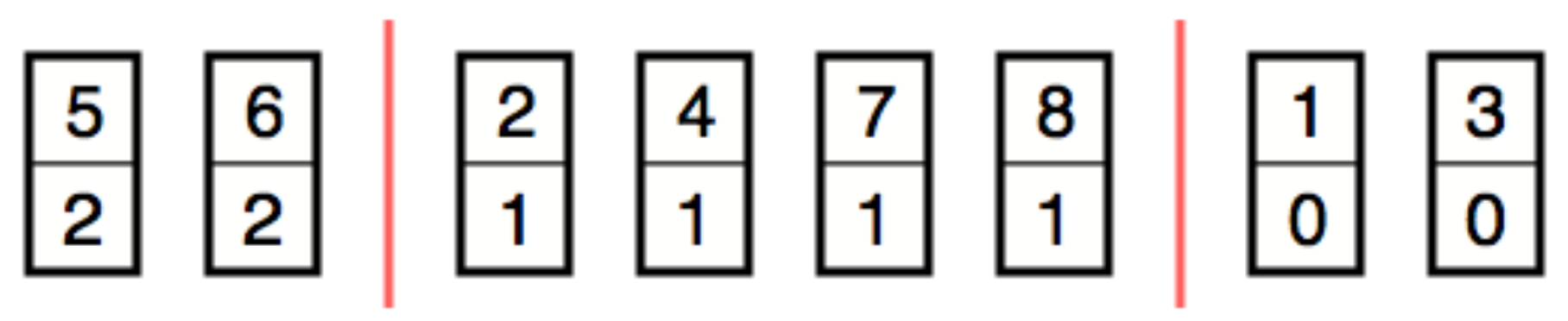}
\vskip -.6in{dominoes($\sig$)}
\end{figure}
\vskip -.14in
\hfill $ 4.4$
\vskip .2in
\noindent
The algorithm constructs all rearrangements of these dominoes that  yield a Parking Function precisely as the array of dominoes on the right of 4.3 may by viewed as yielding the corresponding Parking Function.
For such an arrangement to yield a Parking Function we need to obey the following rules, in placing a domino $[{c_2\atop d_2}]$ immediately to the right of a domino $[{c_1\atop d_1}]$
\begin{align*}
&  \hbox {\it $\bullet$  We must have $d_2\le d_1+1$}\\
&  \hbox {\it $\bullet$  If $d_2=d_1+1$  then we  must also have $c_2>c_1$.}\\
\end{align*}
\vskip -.55in
\hfill $ 4.5$

\vskip .13in
\noindent
The algorithm applied to the domino sequence in 4.4, starts with one of the two  arrangements
\vskip -.08in
{\begin{figure}[H]
\includegraphics[width=.9in]{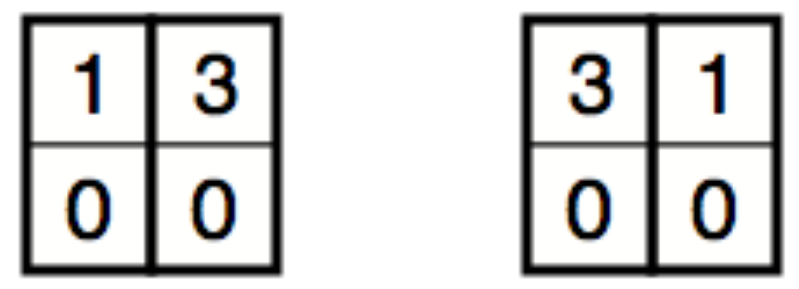}
\end{figure}

\vskip -.43in
\hfill
$4.6$
\vskip .13in
\noindent
and then inserts the rest of the dominoes in 4.4, in the right to left  order, according to 4.5.
\sas

Recall that two cars cause a ``primary dinv'' if they are in the same diagonal and the one on the left is smaller than the one on the right. Alternatively,  two cars cause a ``secondary dinv'' if  the car
on the left is in the immediate higher diagonal than the car on the right and it is greater than the car on the right. To obtain a precise description of the  algorithm let us say that car $a$ {``acts''} on car $b$ if car $b$ is on the right of car $a$ and in the same run,  or car $b$ is in  the run immediately to the right of car $a$ and $b$  is smaller than  $a$. We can easily see that as soon as the domino of car $a$ is placed on the left of a domino of a car $b$ upon which  $a$ acts a new unit of dinv is created.

\newpage

  In the display on the right we have placed
for each car the list of the cars it acts upon. This given, we proceed  as follows, starting from one of the two pairs in 4.6 we construct a tree

\vskip -.065in
\begin{wrapfigure}[8]{R}{.8in}
\vspace{-7pt}
\includegraphics[width=.7in]{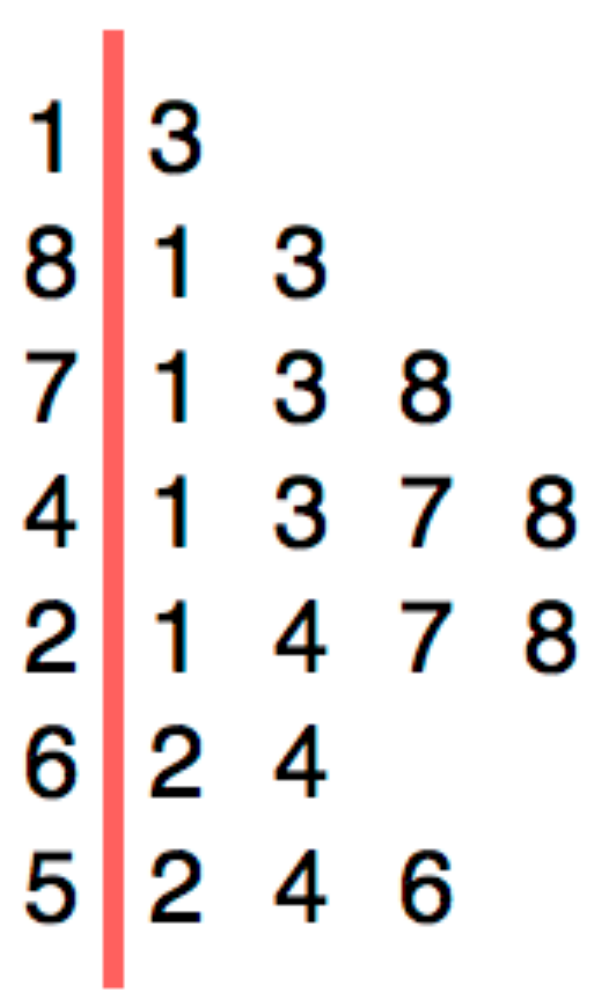}
\end{wrapfigure}

\vskip -.03in
\noindent

 \vskip -.08in

 \noindent
$\CT(\sig)$, whose nodes are indexed by sequences of dominos.  Recursively after the dominoes of cars
$a_1,a_2,\ldots ,a_{i-1}$ have been inserted to obtain a node of $\CT(\sig)$, the labels of the children of that node are obtained by successively inserting in the label of the parent, the domino of car $a_i$ immediately to the right of the domino of one of the cars upon which $a_i$ acts. If this is done proceeding from right to left the additional dinvs that are created are $0,1,\ldots, w_i-1$
where $w_i$ is the number  of cars upon which
$a_i$ acts. For simplicity we will omit the fact that in this  case  car $1$ acts on $3$. In fact, given that we are going to consider only permutations with  last run of length
$2$, car $\sig_{n-1}$ will always act only on $\sig_n$. This given, here and in the following we will denote by $w_i$ the number of
cars acted upon by $\sig_{n-i-1}$. To be consistent with the notation we introduced in the first section, the resulting vector will be denoted $W=(w_1,w_2,\ldots, w_{k-1})$, with $k=n-1$, and will be called a ``schedule'' as before.
In the display below we have depicted an instance of this construction. We have also added on the label of the parent node, a green line to indicate each of the positions the domino of car $2$ can be placed. Finally each   new branch of the tree has been labelled by the power of $q$ corresponding to the resulting dinv increase.
\vskip -.13in
{\begin{figure}[H]
\includegraphics[width=3.4in]{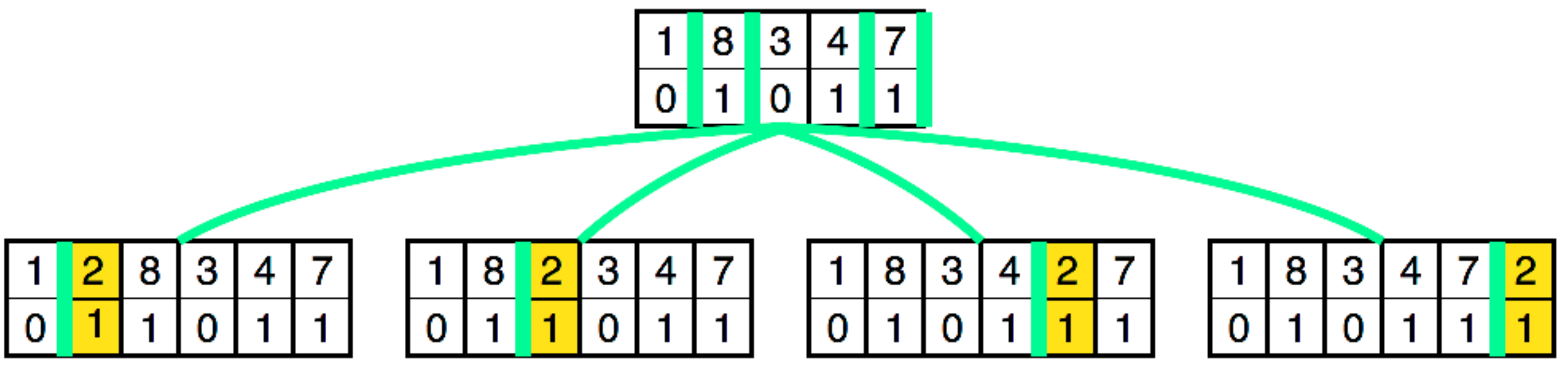}
\end{figure}
\vskip -.17in
The tree obtained when we start from the pair on the left in 4.6  will be called the {\ita left subtree}.
Likewise the tree obtained from the pair on the right of 4.6 will be called the {\ita right subtree}. The leaves of $\CT(\sig)$ yield the domino sequence of the desired Parking Functions
In the following display we illustrate the  choices that yield the leaf  that gives the Parking Function in 4.3
\vskip -.15in
{\begin{figure}[H]
\includegraphics[width=5.2in]{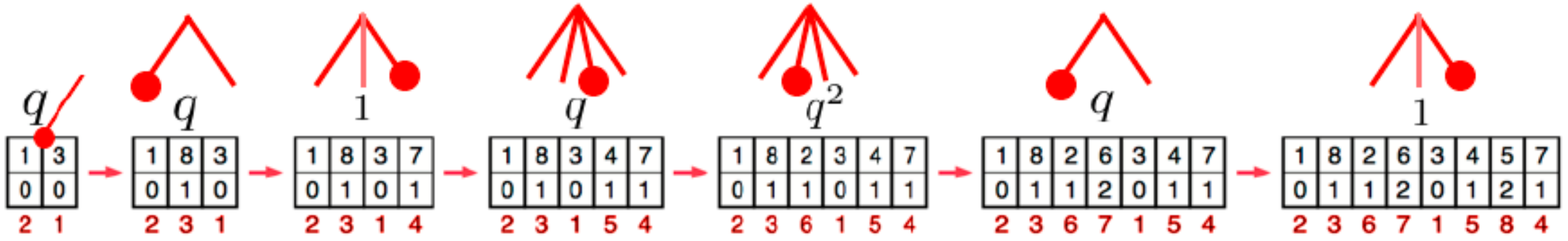}
\end{figure}

\vskip -.2in
The actual path on the  tree of $\sig$ is obtained by appending each step in  the above display right under the red disk of the previous step. The final power of $q$ that gives the dinv of the Parking Function corresponding to   the leaf at the end of the path is simply obtained by multiplying the labels of the branches encountered along the path. To see that the leaves of the resulting tree give all the domino sequences of the Parking Functions that have the desired diagonal car distribution it suffices to see that the algorithm can be reversed by starting from any of the domino sequence of a  desired Parking Functions and removing dominoes of cars
$\sig_1,\sig_2,\ldots ,\sig_n$, as indicated in the following display
\vskip -.1in
{\begin{figure}[H]
\includegraphics[width=5.4in]{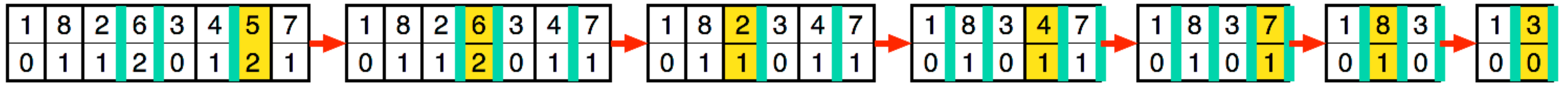}
\end{figure}

\vskip -.18in
It can be easily seen from this that every time we remove a domino its car is necessarily adjacent either to a smaller car in a following run, or a larger car in the same run. It is also  easily noticed  that all the parking functions produced by the tree have the same area statistic. This is due to the fact that every car in diagonal $i$ contributes $i$ area units. Using this information it follows that the common area value is none other than the major index of the permutation $\sig$ that yielded the tree.
This given, for a given $\sig\in  S_n$ whose last run is of size $2$, the  weight of each  $PF\in \CT(\sig)$
 reduces to $q$ to its dinv, the corresponding Gessel fundamental and a power of $x$ to record the diagonal hits
 of its Dyck path. In summary, we are led to define
$$
P_\sig(x;q,F)\ses\sum_{PF\in \CT(\sig)} x^{r(PF)}q^{dinv(PF)}
F_{ides(PF)},
\eqno 4.7
$$
where  $r(PF)$ gives the number of cars weakly to the right of the second diagonal car. In particular the diagonal composition of $PF$  is simply given by the pair  $\big(n-r(PF),r(PF)\big)$.
\newpage

It should be quite apparent now how
Parking Functions can be visually represented by our  Labelled Bar Diagrams.  For instance the Parking Function in 4.3 is  represented by the   labelled bar diagram  in the display below (the last one). The bijection between Parking Functions and
labelled bar diagrams can be easily understood by working on this example.
\vskip -.1in
{\begin{figure}[H]
\includegraphics[width=5.8in]{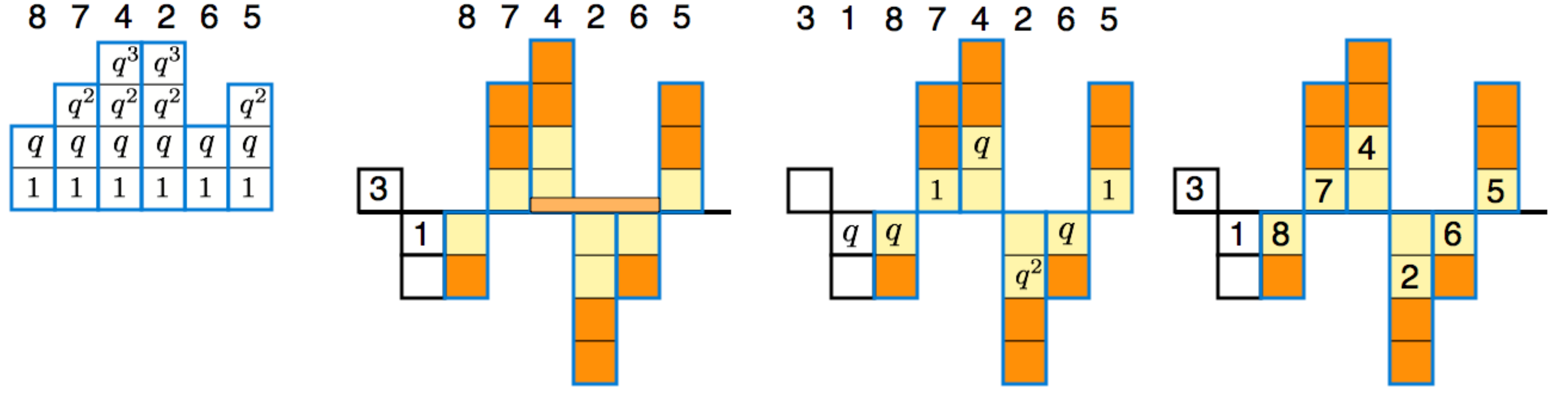}
\end{figure}

\vskip -.1in
\noindent
To begin the length of the bar under car $a$ is simply given by the number of cars acted upon by $a$. The cars are ordered as they appear from right to left in the initial permutation. The cars in the last run are not included in the first diagram in the above display.
The cells are assigned weight as indicated in the table on the left. The second   diagram  in the display is constructed by placing a bar above or below the ground line according as the corresponding car occurs after or before the second car on the main  diagonal. The position of the two cars in the main diagonal is represented by prepending the diagram with
$\UA^3\DA^1 $ or $\DA^3\UA^1 $
according to whether the smaller car  
 precedes or not the larger car.

The reader by now must have at least a glimpse of the understanding how Labelled  Bar Diagrams were created.
In particular we should at least understand that
the powers of $q$ are simply to indicate the amount of dinv  that was created when the corresponding dominos were inserted in the tree $\CT(\sig)$.
More precisely,  if the insertion of the domino of car $a$
created $i$ dinv units then we placed $q^i$ in the cell of that weight in bar $a$. The last diagram is simply obtained by lowering the cars to replace  the powers of $q$. The latter is the labelled bar diagram that   represents the Parking Function in 4.3.

We terminate this section with a list of important properties that easily follow from our construction.
\begin{enumerate}
\item  {\ita The  cars whose bars are below the ground line, in a given bar diagram, precede  the cars whose bars are above  the ground line.}
\item   {\ita A labelled bar diagram   occurs in the construction of all the Parking Functions with prescribed
diagonal cars if and only if none of its bars are totally red.}
\item  {\ita Given a fixed colored unlabeled diagram,
its contribution to the polynomial in 4.7, when we set  $F_{ides(PF)}=1$, is   the polynomial in $q$ obtained by summing all monomials obtained by the
$q$ labelling process described above}
\item {\ita Setting $F_{ides(PF)}=1$ in 4.7 we obtain
 $$
 P_\sig(x;q,1)\ses\sum_{PF\in \CT(\sig)} x^{r(PF)}q^{dinv(PF)}
\ses
\sum_{b=1}^{n-1} A_b(q)\, x^b\ses Q_W(x;q).
\eqno 4.8
$$
Where $n$ gives  the number of cars, and $A_b(q)$ is
the sum of the {\ita dinv} contributions of the  Labelled Bar Diagrams with $b$ bars above the ground line. Here the components of $W$ are given by the successive numbers of cars acted upon by $8,7,2,6,5$.}
\end{enumerate}
\sas

\noindent
Recall that the  operators $C_a$ are defined by setting for any symmetric function $F[X]$
$$
C_aF[X]\ses ( -\tttt{1\over q})^{a-1}F[X\sms \tttt{1-1/q\over z}]\sum_{k\ge 0}z^kh_k[X]\Big|_{z^a}.
\eqno 4.9
$$
It is shown in \cite{HMZ} that  these operators satisfy the relations
$$
q(C_bC_a+C_{a-1}C_{b+1})\ses
C_aC_b + C_{b+1}C_{a-1}
\bigsp \hbox{(for all $b\le a-1$)
}
\eqno 4.10
$$
The compositional Shuffle conjecture by J. Haglund, J. Morse, M. Zabrocki \cite{HMZ}
(now a theorem) states that for any composition $p=(p_1,p_2,\ldots,p_k)$ we have
$$
\nabla C_{p_1}C_{p_2}\cdots C_{p_k}\, 1
\ses
\sum_{{\mathcal PF}(p)}t^{area(PF)}q^{dinv(PF)} F_{ides(PF)}
\eqno 4.11
$$
Where the sum is over all Parking Functions whose Dyck path hits the $0$-diagonal according to $p$.

 The point of departure of the approach in \cite{thesis} was to split the proof of 4.11 by first reducing 4.11 to the case when $p$ is a partition  by means of the identity in 4.10. Then solve the partition case using the fact that in the partition case both sides of  4.11 are symmetric function bases. Leaving the latter case aside for a while let us first see how this reduction
can be accomplished.

To begin notice that
using 4.10 we derive the equality
$$
\nabla C_\aa C_bC_a C_\bb \, 1
\ses  \tttt{1\over q}\nabla C_\aa C_aC_b C_\bb \, 1
\sps \tttt{1\over q}\nabla C_\aa C_{b+1} C_{a-1}C_\bb \, 1\sms
\nabla C_\aa C_{a-1} C_{b+1}C_\bb\, 1
\eqno 4.12
$$
for all $b\le a-1$ and all pairs of compositions $\aa$ and $\bbb$. Now notice that every one of the compositions
$\aa\, a\, b\, \bb$,
$\aa\, b+1\, a-1\, \bb$ and
$\aa\, a-1\, b+1\, \bb$, is lexicographically greater than
$ \aa\,  b\,  a\, \bb$ for $b<a-1$. When $b=a-1$ the resulting identity can still be used to reverse the pair $b<a$. This shows that by a sequence of such uses of the identities
in 4.10 we can express any one of the polynomials $\nabla C_{p_1}C_{p_2}\cdots C_{p_k}\, 1$ as a linear combination of the family of polynomials $\big\{\nabla C_{\la_1}C_{\la_2}\cdots C_{\la_k}\, 1\big\}_{\la\part|p|}$ and since this family is independent, the coefficients are uniquely
determined by the composition $p=(p_1,p_2,\ldots,p_k)$. This given, denoting by $\Pi[p_1,p_2,\ldots,p_k]$ the right hand side of 4.11, we can reduce the proof of 4.11 to the partition case by showing the identities
$$
\Pi[\aa\,  b\,  a\, \bb ]\ses
\tttt{1\over q}
\Pi[\aa\, a\, b\, \bb]\sps
\tttt{1\over q}
\Pi[\aa\, b+1\, a-1\, \bb]\sms
\Pi[\aa\, a-1\, b+1\, \bb]
\eqno 4.13
$$
for all $b\le a-1$ and all pairs of compositions $\aa$ and $\bbb$.

Now it is not difficult to see that 2.5 would be a consequence of the simpler identity
$$
q\big(\Pi[   b\,  a ]+\Pi[  a-1\, b+1]\big)\ses
\Pi[   a\,  b ]+\Pi[  b+1\, a-1]
\eqno 4.14
$$
if the latter could be proved by a bijection of the Parking Functions
contributing to the right hand side
onto those  contributing to the left  hand side that moves  cars only within their diagonals. The reason is very simple, such a bijection would be transferable to a proof of 4.13 since the transfer  could not cause
area changes nor cause changes in dinvs  created by pairs of cars when at least one of the cars is located in the segments covered by $\aaa$ or $\bbb$. For the same reason  the $ides$ of a $PF$ could not be affected under the transfer again for any pair of  cars $i,i+1$ when at least one of them is in the segments covered by $\aa$ or $\bb$.
\sas

The miracle here is the further discovery that we can prove that such a bijection 
 already exits within the collections of Parking Functions obtained from any permutation
$\sig $ with a last run of size $2$.
More precisely, given $\sig\in  S_n$
with a last run of length $2$, if we set
$$
\Phi_\sig[a,b]\ses
\sum_{ PF\in  {\mathcal T}(\sig)}q^{dinv(PF)}
F_{ides(PF)}\chi\big(p(PF)=(a,b)\big),
$$
we can prove  the validity of the identity
$$
q\big(\Phi_\sig[b,a]\sps \Phi_\sig[a-1,b+1]\big)= \Phi_\sig[a,b]\sps \Phi_\sig[b+1,a-1]
\bigsp
\hbox{
(for   $b\ge 1$  and  $ a\ge 1$
)}
\eqno 4.15
$$
This given, in principle, the desired reduction could be achieved by constructing, for any given such $\sig \in S_n$, a bijection of the family Parking Functions contributing to the left hand side of 4.15 onto the family contributing to the right hand side which preserves   {\ita ides} and
increases   {\ita dinv} exactly by one.
Now it is easily seen that for such $\sig$ the polynomial in 4.7 has the expansion
\vskip -.18in
$$
P_\sig(x;q,F)=\sum_{PF\in \CT(\sig)} x^{r(PF)}q^{dinv(PF)}
F_{ides(PF)}
\ses\sum_{k=1}^{n-1}x^k A_k(q,F)
\eqno 4.16
$$

\noindent
Thus, making the replacements
$a=s+1$, $b=n-a$ and $n=k+1$,   4.15,
can  also rewritten as
$$
q\big(A_{s+1}(q,F)+A_{k -s +1}(q,F)\big)
=A_{k-s}(q,F)+A_{s}(q,F)
\bigsp
\hbox{
(for   $  k-s\ge 1$ and   $s\ge 1$
)}
\eqno 4.17
$$
At  this point the existence of the desired bijection seems to require a proof that  the  coefficients of all the polynomials  $P_\sig(x;q,F)$  satisfy the identities in 4.17.

\newpage

That may appear as a tall order.  Fortunately, another rather surprising miracle comes to our help. Guided by computer data we discovered that is relatively easy to establish in full generality a factorization of the form
$$
P_\sig(x;q,F)=P_\sig^{(1)}(x;q)\times
P_\sig^{(2)}(q;F)
\eqno 4.18
$$
To understand what causes   such a factorization,
we need some definitions. A string of consecutive integers $i,i+1,\cdots,i+\l-1$ in a given run with both $i-1$ and $i+\l$ not in this run, is  called a {\ita maxicon}.
An index $i$ that is in the $ides(PF)$ of every leaf of $\TAU(\sig)$ is  called a {\ita forced
ides}. The following observations are immediate consequences of our construction of the polynomials  $P_\sig(x;q,F)$:
\vskip -.2in

\begin{enumerate}
\item{\ita The elements of a maxicon  $i,i+1,\ldots ,i+\l-1$  will appear in the leaves of
$\TAU(\sig)$ in each of their $\l!$ orders}
\item{\ita If a given car $i$ belongs to a run that follows the run that contains $i+1$, then $i$ will always be in the $ides(PF)$ of every leaf of    $\TAU(\sig)$  }
\end{enumerate}
This given, the factorization in 4.18, proved in \cite{thesis}, may be stated as follows.
\sas

\noindent{\bol Theorem 4.1}

{\ita Calling ``Ycons`` the Young subgroup of $S_n$ generated by the maxicons of a $\sig \in S_n$ with a last run of length $2$, we have the factorization
$$
P_\sig(x;q,F)\, =\,
\Bigg(
\sum_{\multi {PF\in \CT(\sig) \cr ides(PF)=fides(\sig)}}  x^{r(PF)}q^{dinv(PF)}
\Bigg)
\bigg(\sum_{\aaa\in {Y\hskip -.03in cons}}q^{inv(\aaa)}F_{ides(\aaa)\cup fides(\sig)}
\bigg)
\eqno 4.19
$$
where for convenience we have let $fides(\sig)$ denote the set of forced ides of $\CT(\sig)$.}
\sas

 This is all beautifully confirmed by the following display that constructs  the tree $\TAU(\sig)$ for  $\sig=45312$. The runs of $\sig$ are $4\,5$, $3$ and $1\, 2$. Thus we have  $2$ maxicons
 $1,2$ and $4,5$.  Since $2$ is in the main diagonal, $3$ is in diagonal $1$ and
$4$ is in diagonal $2$, we deduce that in the word of each leaf of $\TAU(\sig)$ we will have $4$ to the left of $3$ and $3$ to the left of $2$. Thus $3$ and $2$ will be both forced ides. In this case the tree yields
$$
P_\sig(x;q,F)=
x^4(q^3+q^2)F_{1,2,3,4}\sps
x^4(q^2+q)F_{1,2,3 }\sps
x(q^2+q)F_{ 2,3,4}\sps
xF_{2,3}
$$
whose factorization is at the end of the display
\vskip -.03in
{\begin{figure}[H]
\includegraphics[width=5.3in]{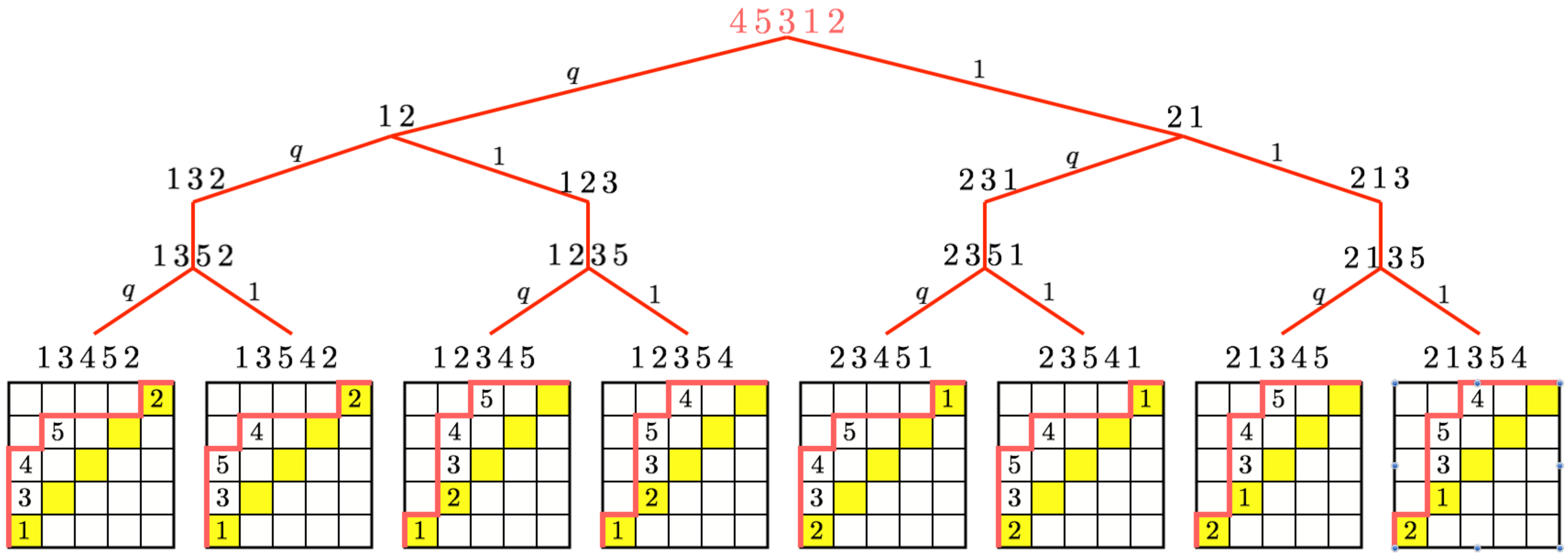}
\end{figure}
\vskip -.28in
{\begin{figure}[H]
\includegraphics[width=5.3in]  {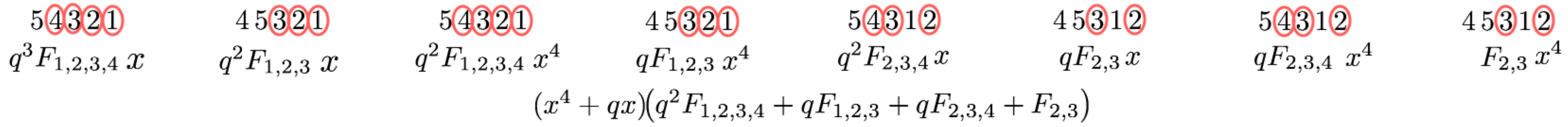}
\end{figure}
\vskip -.1in

In particular the validity of the factorization in 4.18
shows that 4.17 may be rewritten as
$$
q\Big(
P_\sig^{(1)}(x;q)\big|_{x^{s+1}}
+
P_\sig^{(1)}(x;q)\big|_{x^{k-s+1}}
\Big) P_\sig^{(2)}(q;F)
-\Big(
P_\sig^{(1)}(x;q)\big|_{x^{k-s}}
+
P_\sig^{(1)}(x;q)\big|_{x^s}
\Big)P_\sig^{(2)}(q;F).
\eqno 4.20
$$
Rather than canceling the common factor
we make the replacement $F\RA 1$
and, using 4.18,  rewrite 4.19 as
$$
q\Big(
P_\sig  (x;q,1)\big|_{x^{s+1}}
+
P_\sig (x;q,1)\big|_{x^{k-s+1}}
\Big)
=
P_\sig (x;q,1)\big|_{x^{k-s}}
+
P_\sig(x;q,1)\big|_{x^{s}}
$$
\newpage

Now it is not difficult to see that, in full generality,  we have the equality
$$
P_\sig (x;q,1)\ses Q_W(x,q),
\eqno 4.21
$$
where for $\sig=\sig_1\sig_2\cdots \sig_n$ a permutation with last run of length $2$
the components of $W$ are obtained by the now quite familiar algorithm. That is we first break up  $\sig$ into runs set
$x_i=\sig_{n-i-1}$, and, if  $i\ge 1$, set $w_i$ to be the number of $x_j>x_i$   that are either in the same run as $x_i$ or $x_j<x_i$ and in the following run.
Notice further that, with  this convention, $x_0=\sig_{n-1}<\sig_n=x_{-1}$. Thus all the Parking Functions with car $x_{-1}$  to the left of car $x_0$,
will necessarily have and extra primary dinv caused by these two cars. That explains the extra factor $q$ we
assigned to all the monomials produced by the Labelled Bar Diagrams on the left subtree of $\TAU(\sig)$. It is easy to see that this creates schedules which can increase by at most one, and that all legal schedules can be constructed this way. 

In summary, the bijection between  Parking Functions yielded by  the leaves of $\TAU(\sig)$ and Labelled
Bar Diagrams should by itself make evident the equality in 4.21. Moreover, our proof of the functional equation, combined with the identity in 4.21, yields the identity
 \vskip -.2in
 $$
 P_\sig(x;q,1)\ses\sum_{PF\in \CT(\sig)} x^{r(PF)}q^{dinv(PF)}
\ses Q_W(x;q)\ses
\sum_{b=1}^{n-1} A_b(q)\, x^b,
\eqno 4.22
$$

\vskip -.2in
\noindent
with
$$
q\big(A_{s+1}(q )+A_{k -s +1}(q )\big)
=A_{k-s}(q)+A_{s}(q),
\bigsp
\hbox{
(for   $  k-s\ge 1$ and   $s\ge 1$.
)}
\eqno 4.23
$$
Thus by combining 4.21 with the proof of the functional equation we derive the following fundamental fact.
\sas

\noindent{\bol Theorem 4.2}

{\ita The equality in 4.11, that is
$$
\nabla C_{p_1}C_{p_2}\cdots C_{p_k}\, 1
\ses
\sum_{{\mathcal PF}(p)}t^{area(PF)}q^{dinv(PF)} F_{ides(PF)}
\eqno 4.24
$$
holds for all compositions if and only if it holds for all partitions}
\sas

In other words, to complete the proof of the Compositional Shuffle Theorem we need to show the equality
$$
\nabla
C_\la\, 1[X;q,t]\ses
\Pi_\la[X;q,t]
\eqno 4.25
$$
for all partitions $\la$.
\sas

It is easy to derive that the family of polynomials on the left hand side of 4.25 are a symmetric function basis. The same can be proved also for the family on the right hand side. This given, our present plan is to derive the equality in 4.25 as an application of the  following  elementary linear algebra result. In fact, this theorem   turned out to be successful  in the proof  of the Haglund combinatorial formula for the modified Macdonald Polynomials \cite{HagHaiLo} as well as for the balanced path result in \cite{Hag08}.
\sas

\noindent{\bol Theorem 4.3}

{\ita Suppose an $n$ dimensional  vector space $\BV$ has two bases
$$
\LL \phi\RR=\LL \phi_1,\phi_2,\ldots,\phi_n\RR,
\bigsp\bigsp
\LL \psi\RR=\LL \psi_1,\psi_2,\ldots,\psi_n\RR.
$$
This given, any basis
 $\ess
\LL H\RR= \LL H_1,H_2,\ldots,H_n\RR
\ess $
which is Upper triangularly related to the $\LL\phi\RR$ basis and Lower triangularly related to the $\LL\psi\RR$ basis is uniquely determined by the normalizing condition
$$
F[H_i]= c_i \neq 0
\bigsp (\hbox{  for $1\le i\le n$ })
\eqno 4.26
$$
for a suitable linear functional $F$}
\sas

It turns out that that in \cite{GXZ} it is shown that
for all $\la\part n $ we have
$$
\LL \nabla
C_\la\, 1[X;q,t]\scs e_n\RR\ses
\LL\Pi_\la[X;q,t]\scs e_n\RR.
\eqno 4.27
$$
Thus we can use this result to provide the desired functional $F$ in 4.26. It is not difficult to prove the additional triangularity results for the symmetric function  side with
$$
a)\ess\ess\LL \phi\RR\,=\, \LL s_\la[\tttt{X\over q-1}]\RR_{\la \part n}\bigsp
 \hbox{and}\bigsp
b)\ess\ess\LL \psi\RR\,=\, \LL \nabla s_\la[X]\RR_{\la \part n}
\eqno 4.28
$$
\newpage

However, so far we have only been able to prove that the combinatorial side is upper triangularly related to the basis in  4.28 a). But the lower triangularity of the combinatorial side to the basis in 4.28 b) is
unlikely to lead to a proof since the still conjectural \cite{LoWar} combinatorial interpretation of
$\nabla s_\la$  has $\nabla e_n$ as a special case.

The upper triangularity result,
for both bases, which is highly non trivial for the combinatorial side, will be presented in a forthcoming publication.
\sas

In the next few pages we will first derive some further consequences of the functional equation  and formulate an interesting open problem.
We will terminate the paper with some conjectures that  sharpen the results of this paper.
\sas

Notice first that the factorization in 4.19  enables us to recover the second factor in 4.18 and thus also the whole polynomial $P_{\sig}(x;q,F)$.
 But it does not stop here. Indeed, there is no reason to limit our Labelled Bar Diagrams to representing leaves of $\TAU(\sig)$ for $\sig$ a permutation with last run of length $2$. In fact, we need only make a few  changes when we place no restriction on the length of the last run. More precisely, for a $\sig \in S_n$, a permutation with last run of length $r$
the corresponding, polynomials  can be  written in the form
$$
Q_\sig(Y_r;q,F)\ses \bu  \sum_{p=(p_1,p_2, \ldots, p_r)\models n}\bu
y_1^{p_1}y_2^{p_2}\cdots y_r^{p_r}\, \Pi_\sig(p; q,F)
\eqno 4.24
$$
with $Y_r=(y_1,y_2,\ldots , y_{r})$, where the sum in 4.24 is over $r$-part compositions of $n$;

 For such $\sig$ the only  difference in the construction of the  tree $\TAU(\sig)$, is that the root of this tree will have $r!$ children, according to the
permutation of the cars in the last run. It goes without saying that the $r!$ branches of $\TAU(\sig)$ emanating from its root are labelled with powers of q summing to the polynomial $[r]_q!$. Of course, the Parking Functions yielded by the leaves of $\TAU(\sig)$ will all have a  diagonal composition of length $r$.
In particular, we have  the expansion
$$
\Pi_\sig(p; q,F) \, = \bu \sum_{PF\in {\mathcal PF}_{\sig}(p)}\bu q^{dinv(PF)}F_{ides(PF)}
\eqno 4.25
$$
where the sum is over all Parking Function with diagonal cars given by the runs of $\sig$ and whose
Dyck path hits the main diagonal according to the composition $p$.
\sas

Now suppose that the composition $p$ may be decomposed in the form $p=\aa\, b\,a\, \bb$ with $\aa$ and $\bb$ such  that $p\models n$ but otherwise arbitrary, and $1\le b\le a-1$. Then it follows from the functional equation that we must have the identity
\begin{align*}
q\big(
\Pi_\sig(\aa\, b\,a\, \bb\, ;\, q,F)&
\sps
\Pi_\sig(\aa\, (a-1)\,(b+1)\, \bb\,;\, q,F)
\big)
\ses
\\
&\ses
\Pi_\sig(\aa\, b\,a\, \bb\, ;\,  q,F)
\sps
\Pi_\sig(\aa\, (a-1)\,(b+1)\, \bb\, ;\, q,F)
\end{align*}
\vskip -.3in
\hfill
4.26

\vskip .1in
\noindent
and that must holds true for every $\sig\in S_n$
and for all such  $\aa,\bb$ and $1\le b\le a-1$.
\sas

The presence of all these identities, even when the Parking Functions  are restricted to have prescribed diagonal cars, strongly suggests that all the ingredients occurring in the Symmetric Function
world must have corresponding specializations to the Quasisymmetric Funtions world.
\sas

To be precise we are led to conjecture  that there must exist Quasisymmetric analogs of Modified Macdonald polynomials, nabla as well as the $C_a$ operators satisfying the fundamental identity in 4.10. Recall that in the Symmetrc Function world all the  composition  identities follow from 4.10.
We find it difficult to believe that the phenomenon
resulting from the functional equation has no Quasisymmetric function counterpart. Several attempts have been made in the past to carry
out such extensions. The first one that seems to  be worth   investigating  is the one given by N. Bergeron and M. Zabrocki in  \cite{MikeNant}.
\sas

But there are even more amazing facts that strengthen the validity of these conjectures. Our point of departure is here is the following identity:
\newpage

\noindent{\bol Proposition 4.1}

{\ita For all $a>b+1$ we have}
$$
(-\tttt{1\over q})^{a+b-3} s_{a-1,b+1}[X]
\ses \BC_a\BC_b\, {\bf 1}\sms q \BC_{a-1}\BC_{b+1}\, {\bf 1}
\eqno 4.27
$$
\noindent{\bol Proof}

From 4.9 it follows that
$$
C_b\, 1\ses (-{\TS{ 1\over q }})^{b-1} h_b[X]
$$
Using 4.9 again we get
\begin{align*}
(-q)^{a+b-2}C_a C_b\, 1
&\ses
 h_b[X+\tttt{1-  q \over q z}]
\sum_{k\ge 0}z^k h_k[X]\, \Big|_{z^a}\ses h_b h_a \sps
(1-q)\sum_{r\ge 1}h_{b-r}h_{a+r}/q^r
\end{align*}
Likewise we obtain
\begin{align*}
(-q)^{a+b-2}qC_{a-1} C_{b+1}\, 1
&\ses qh_{b+1} h_{a-1} \sps
(1-q)\sum_{r\ge 1}h_{b+1-r}h_{a-1+r}/q^{r-1}
\end{align*}
Thus
$$
(-q)^{a+b-2}
(C_a C_b\, 1-qC_{a-1} C_{b+1}\, 1)
=
h_b h_a -qh_{b+1} h_{a-1} \sms (1-q)
h_bh_a= q(h_b h_a-h_{b+1} h_{a-1} )
$$
and the Jacobi-Trudi identity gives 4.27.
\sas

Now combining 4.27 with the conjectures in \cite{BGHT} for Nabla of a Schur function  we derive the following
\sas

\noindent{\bol  Conjecture III}

{\ita The symmetric polynomial
$$
\nabla  \BC_a\BC_b\, {\bf 1}\sms q \nabla \BC_{a-1}\BC_{b+1}\, {\bf 1}\ses
q\nabla  \BC_b\BC_a\, {\bf 1}\sms   \nabla \BC_{b+1}\BC_{a-1}\, {\bf 1}
\eqno 4.28
$$
is Schur positive.
}
\sa

Notice that applying Nabla to both sides of 4.27 we can rewrite it in the form
$$
\nabla s_{a,b}\ses (-q)^{a+b-3}\Big(
\nabla C_{a+1}C_{b-1}{\bf 1}\sms q\nabla C_aC_b{\bf 1}
\Big)
$$

This given, we should mention that Yeon Kim 
in \cite{yeon}, in both cases $(a,b)=(n-3,3)$ and 
$(a,b)=(n-4,4)$, succeeded to construct  an injection of the Parking Functions with diagonal composition $(a,b)$ 
into Parking Functions with diagonal composition 
$(a+1,b-1)$ preserving area and Gessel Fundamental 
and increasing dinv by one unit. Since  her work 
$n$ was only required to be greater than $4$ and $7$ respectively in the above two cases she succeeded in giving a Parking Function interpretation to the resulting infinite variety 
of Nabla Schurs $\nabla s_{n-3,3}$ and $\nabla s_{n-4,4}$.
\sas

In this connection,  computer experimentation 
revealed  that Kim's type results should be obtainable even with prescribed  diagonal cars, \

For a more precise description of these findings  we must recall the polynomial
$$
\Phi_\sig[a,b]\ses
\sum_{ PF\in  {\mathcal T}(\sig)}q^{dinv(PF)}
F_{ides(PF)}\chi\big(p(PF)=(a,b)\big).
$$
It turns out that computer data strongly suggests  that the identity in 4.15 can be sharpened to
$$
\Phi_\sig[a,b]\sms q\, \Phi_\sig[a-1,b+1]\ses
q\, \Phi_\sig[b,a]\sms \Phi_\sig[b+1,a-1]\in \BN[q]
\eqno 4.29
$$
for all $\sig\in S_n$, $a+b=n$ and  $a,b\ge 1$.
\sas

\noindent
From 4.29 it follows that for any composition $\ggg\models d $ we have
$$
\Big(\Phi_\sig[a,b]\sms q\, \Phi_\sig[a-1,b+1]\Big)\Big|_{F_\ggg}\ses
\Big(q\, \Phi_\sig[b,a]\sms \Phi_\sig[b+1,a-1]\Big)\Big|_{F_\ggg}\in \BN[q]
\eqno 4.30
$$
This is a remarkable parallel to
Conjecture III at the Quasi-symmetric function level. It also suggests for instance that there is an injection of Parking Functions with diagonal cars prescribed by $\sig$ and diagonal composition
$a-1,b+1$ into  Parking Functions  with diagonal composition $a,b$ which preserves the Gessel fundamental and increases the dinv by one unit.
Likewise there must be a way of injecting
Parking Functions with diagonal cars prescribed by $\sig$ and diagonal composition
$b+1,a-1$ into  Parking Functions  with diagonal composition $b,a$ which preserves the Gessel fundamental and decreases  the dinv by one unit.
\sas

But there is one more surprising sharpening of the functional equation by passing to our schedule polynomilals $Q_W(x;q)$. To make easier to compare our statements with 4.30, let  us  write this polynomial,
 for the schedule $W=(w_1,w_2,\ldots , w_{k-1})$, in the form
$$
Q_W(x,q)\ses \sum_{\multi{b\in [1,k]\cr a+b=k+1\cr}}
A_{a,b}^W(q)\, x^b
\eqno 4.31
$$
This given, passing from 4.16 to 4.31 using the factorization result, 4.30 becomes
$$
 A^W_{a,b}(q)\sms q\, A^W_{a-1,b+1}(q)
 \ses
 q\, A^W_{b,a}(q)\sms A^W_{b+1,a-1}(q) \in \BN[q]
 \ess\ess\ess \hbox{( for all $1\le b\le k$)}
\eqno 4.32
$$
To see what this implies let us explore the case  $n=7$.
This gives

\vskip -.15in

\begin{align*}
&a)\ess\ess   A^W_{6,1}(q)\sms q\,A^W_{5,2}(q)
\ses q\, A^W_{1,6}(q)\sms A^W_{2,5}(q) \ses q\, A(q)
\in \BN[q],\\
\ess\ess\ess
& b)\ess\ess   A^W_{5,2}(q)\sms q\,A^W_{4,3}(q)  \ses q\, A^W_{2,5}(q)\sms A^W_{3,4}(q)
\in \BN[q],\\
\ess\ess\ess
&c)\ess  q\, A^W_{3,4}(q)\ses A^W_{4,3}(q)
\in \BN[q].
\end{align*}

\vskip -.45in

\hfill  4.33

\vskip .3in
\noindent
Now a) implies that the polynomial $A^W_{2,5}(q)$
must be divisible by $q$. Likewise b) implies that
$A^W_{3,4}(q)$ must be divisible by $q^2$ and finally
c) forces $A^W_{4,3}(q)$  to be divisible by $q^3$.
Thus writing
$$
a)\ess\ess A^W_{2,5}(q)=q\, A^{W,1}_{2,5}(q)
\ess\ess\ess\ess\ess\ess
b)\ess\ess A^W_{3,4}(q)=q^2\,C(q)
\ess\ess\ess\ess\ess\ess
c)\ess\ess A^W_{4,3}(q)=q^3\,C(q)
$$
the identities in  4.33 can be rewritten as
$$
a)\ess\ess q\, A^W_{1,6} - q\, A^{W,1}_{2,5} =q\, A(q)   \ess\ess\ess\ess\ess\ess
b)\ess\ess q^2 A^{W,1}_{2,5} - q^2\,C(q) =q^2\, B(q)
$$
In summary we have
$$
c)\ess\ess A_{3,4}^W\ses q^2\, C
\ess\ess\ess\ess\ess\ess
b)\ess\ess   A_{2,5}^W\ses q\, B \sps q\, C
\ess\ess\ess\ess\ess\ess
a)\ess\ess A_{1,6}^W\ses A\sps B\sps C.
\eqno 4.34
$$
Thus from  4.33 it follows that
$$
a) \ess\ess A_{6,1}^W-q\,A_{5,2}^W=
q\,(A+B+C)- (qB +qC)
\ess\ess\ess\ess\ess\ess
b) \ess\ess A_{5,2}^W-q\,A_{4,3}^W=
q^2\, B+q^2\,C-q^2\,C
$$
and using 4.33 c) we finally obtain
$$
a)\ess\ess A_{6,1}^W=q\, A+q^3\, B+q^5\, C
\ess\ess\ess\ess\ess\ess
b)\ess\ess A_{5,2}^W= q^4\,  C+q^2\, B
\ess\ess\ess\ess\ess\ess
c) \ess\ess A_{4,3}^W= q^3\,  C.
\eqno 4.35
$$

 \vskip -.03in
\begin{wrapfigure}[11]{R}{2.4in}
\vspace{-7pt}
\includegraphics[height=1.65in]{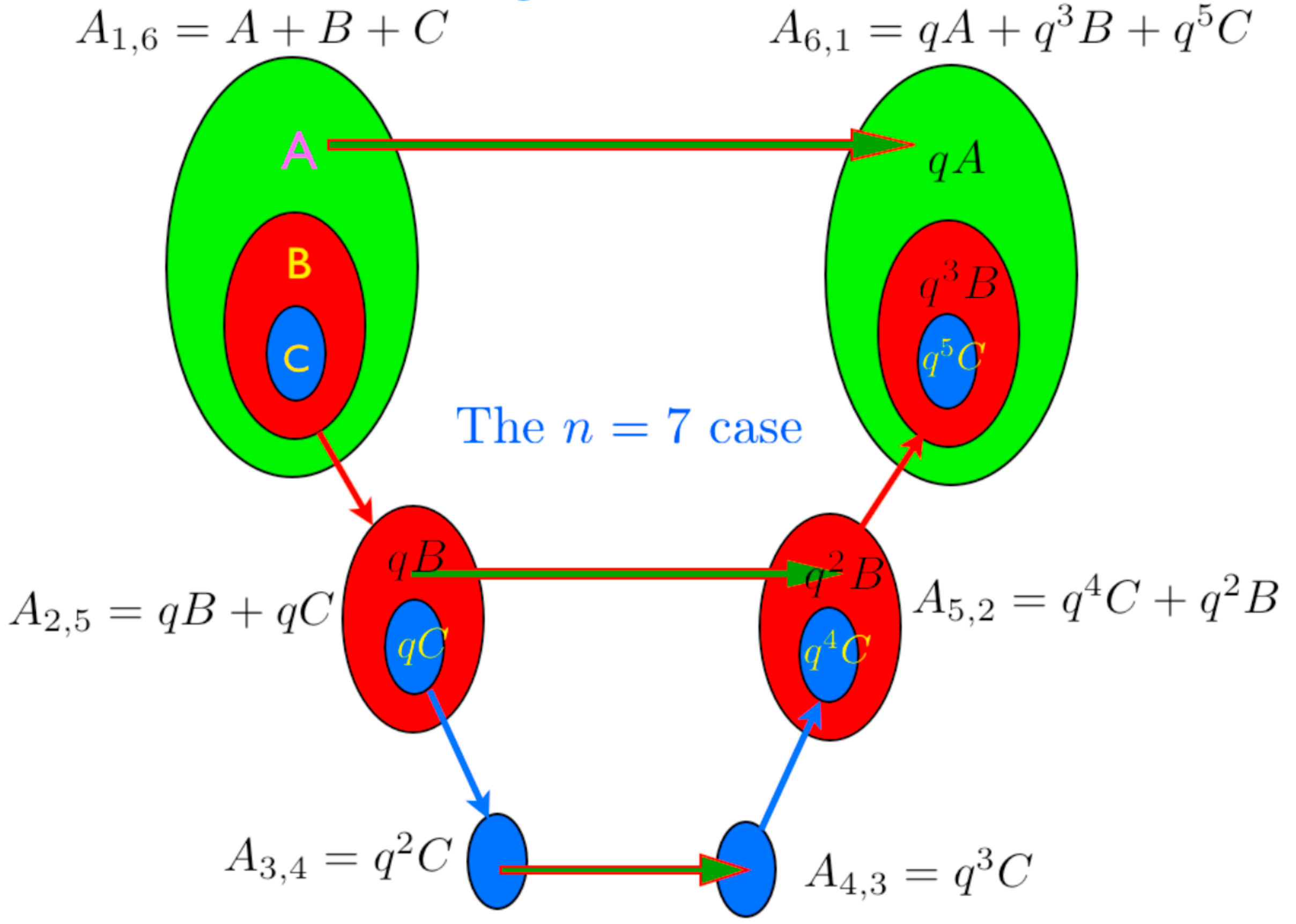}
\end{wrapfigure}

\vskip -.1in
\noindent

\vskip .02in

\noindent
This is typically what happens when $n$ (the number
of cars) is odd. Our data not only confirms
all of these identities but reveals that the polynomials $A,B,C$ are also unimodal. The implications of these identities are best understood  by  a visual display. In the figure on th right we have depicted the polynomials
$A(q),B(q)$ and $C(q)$ respectively in Green, Red and Blue. The Red down arrow represents a conjectural bijection of a sub polynomial  of
$A_{1,6}^W$ onto $A_{2,5}^W$ which increases the dinv by one unit. The Red up arrow on the right
represents a conjectural bijection of $A_{5,2}^W$
onto a sub polynomial  of $A_{6,1}^W$   which increases the dinv by one unit. The top Brown arrow represents a highly non trivial bijection constructed by our second author during her  thesis
work, and played a fundamental role in \cite{AHicks}. It may not be unlikely that all the
intricacies encountered in the proof of this bijection could be explained by means of this display. The middle  Brown arrow  represents a   bijection  of $A_{2,5}^W$ onto $A_{5,2}^W$ whose  proof in full generality appears difficult. The bottom Brown arrow should be easy
to prove.
\newpage

 The following basic result was proved in   lecture notes on schedule polynomials:  
\sa

\noindent{\bol Theorem 4.4}

{\ita The vector space of polynomials
\vskip -.25in
$$
Q(x,q)\ses \sum_{i=1}^k A_i (q)x^i
\eqno 4.36
$$
\vskip -.15in
\noindent
which satisfy the functional equation
$$
(1-q/x)Q(x)\sps x^k (1-qx)Q(1/x) \ses (1+x^k)
\big(A_k(q)-q A_1(q)\big)
$$
is $\lfloor (k+1)/2\rfloor$ dimensional with basis the following polynomials  }

{a)} {\ita When  $k=2a$ }
\vskip -.25in
$$
E_{j,k}(x;q)=
\sum_{i=1}^j\big (q^{2j-i}x^i\sps q^{i-1} x^{2a+1-i}  \big)
\ess\ess\ess\ess\ess\ess
\hbox{for $1\le j\le a$}
$$

{b)}$\hskip .012in$ {\ita When  $ k=2a-1 $ }
\vskip -.2in
\begin{align*}
O_{j,k}(x;q)=
 \sum_{i=1}^{j }\big(q^{2j-i}  x^i
 &\sps    q^{i-1}   x^{2a -i}\big)
\ess\ess\ess\ess \hbox{for $1\le j\le a-1\ess $,}
\end{align*}
\vskip -.2in
\begin{align*}
\ess\ess\ess\bigsp\bigsp\bigsp
\hbox{\ita and } \ess\ess
\ess O_{a,k}(x;q)=\sum_{i=1}^{a-1}\ q  ^{2a-i-1}  x^i\sps  \sum_{i=1}^{a}    q^{i-1}   x^{2a -i}.
\end{align*}

\noindent
Notice, for instance that from a) it follows that

\vskip .1in
 $
E_{1,6}= x^6+q   x, \ess\ess
E_{2,6}= x^6+q  x^5 +q^2  x^2+q^3  x,\ess\ess  $
\vskip -.15in
$$
\bigsp\bigsp\bigsp
E_{3,6}= x^6+q x^5 +q^2  x^4+q^3  x^3 +
q^4  x^2+  q^5x
\eqno 4.37
$$
The independence of the basis in Theorem 4.4 is due to a very simple reason, easily seen in this  example.
 Namely the fact,  that both families
 $\{E_{j,2a}(x;q)\}_ {1\le j\le a}$ and
  $\{O_{j,2a-1}(x;q)\}_ {1\le j\le a}$
are  triangularly related to the monomials $x,x^2,\cdots ,x^a$.
\sas

For any legal schedule $W=(w_1,w_2,\ldots,w_{5})$ the coefficients $A^W_i(q)$ of  the polynomial
$$
Q_W(x;q)\ses \sum_{i=1}^6 A^W_i(q)\, x^i
\eqno 4.38
$$
must satisfy the identities in 4.33. Thus it follows that it may be rewritten in the form
\begin{align*}
Q_W(x;q)=&
(A+B+C)x^6+
(qB+qC)x^5+
q^2 Cx^4+\\
&\bigsp\ess\ess
q^3Cx^3+
(q^2 B+q^4C)x^2+
(qA+q^3B+Cq^5)x.
\end{align*}
This fact is an immediate consequence of 4.34 and 4.35. Regrouping terms according to   $A,B,C$  gives 
\begin{align*}
Q_W(x;q)=&
A(q)\,(x^6+qx) \sps
B(q)\,(x^6+qx^5 +q^2x^2+q^3  x) \sps\\
&\bigsp\bigsp\bigsp
\sps C(q)\,(x^6+q x^5 +q^2  x^4+q^3  x^3 +
q^4  x^2+  q^5x)
\end{align*}

\noindent
Comparing the coefficients of $A,B,C$ with 
 the basis in 4.37 gives the expansion
\begin{align*}
Q_W(x;q)=&
A(q)\,E_{1,6}(x;q) \sps
B(q)\,E_{2,6}(x;q) \sps
C(q)\,E_{3,6}(x;q)
\end{align*}
The computations we carried out
 in this particular case
should clearly indicate how the bases introduced in a) and b) of Theorem 4.4 were discovered.
\sas

The algorithm that constructs the coefficients of the expansion of each schedule polynomials in terms of the appropriate basis proceeds  very much as was illustrated in the special case $k=6$. More precisely, given a schedule
$W=(w_1,w_2,\ldots ,w_{k-1})$, let
$$
Q_W(x;q)\ses\sum_{s=1}^k A_{k+1-s,s}(q)\, x^s
$$

For $k=2a=n-1$  set
$$
\aa_r(q)=
q^{-r}\big(qA_{r,n-r}(q)- A_{r+1,k-r}(q) \big)
\ess\ess\ess\hbox{for  $1\le r\le a-1\ess\ess \&$ }\ess
\aa_a(q)=q^{-a+1}A_{a,a+1}(q)
\eqno 4.39
$$

For $k=2a-1$  set
$$
\aa_r(q)=
q^{-r}\big(qA_{r ,n-r}(q)\sms A_{r+1,k-r}(q) \big)
\ess\ess\ess\hbox{for  $1\le r\le a-1\ess\ess \&$ }\ess
\aa_a(q)=q^{-a+1}A_{a,a}(q)
\eqno 4.40
$$
The factor $q^{-r}$  can be shown to be the  minimum power  of $q$ occurring  in the polynomial
$qQ_{r}(q) $.

This given, we have
\sas

\noindent{\bol Theorem 4.5}

{\ita Using 4.39 or 4.40, as the case may be,
the polynomial $Q_W(x,q)$ has the expansion }
$$
Q_W(x;q)\ses
\begin{cases}
  \sum_{r=1}^a \aa_r(q)\,  E_{r,k}(x;q) &
\text{ if}\ess\ess  k=2a  \cr
\\
  \sum_{r=1}^a \aa_r(q)\,  O_{r,k}(x;q) & \text{if}\ess\ess  k=2a-1
\end{cases}
  \eqno 4.41
$$
\noindent{\bol Proof}

In the first case we need only show that
$$
a)\ess\ess A_{s,n-s}^W=q^{s-1}(\aa_s+\aa_{s+1}+\cdots+ \aa_a)
\ess\ess\ess\ess\ess
\hbox{and}
\ess\ess\ess\ess\ess
b)\ess\ess A_{n-s,s}^W= \sum_{s\le r\le a}q^{2r-s}\aa_r.
$$
While in the second case we need only show

$$
a)\ess\ess A_{s,n-s}^W=q^{s-1}(\aa_s+\aa_{s+1}+\cdots+ \aa_a)
\ess\ess\ess\ess\ess
\hbox{and}
\ess\ess\ess\ess\ess
b)\ess\ess A_{n-s,s}^W=
\sum_{s\le r\le a-1}q^{2r-s}\aa_r\sps q^{2a-s-1}\aa_a .
$$

This can be carried out in a straightforward manner. Since we already illustrated the $k$ even case, we will limit our efforts to verifying  the odd case $k=7$. Using 4.40  for $a=4$  we obtain
\begin{align*}
\aa_1(q)= q^{-1}\big(qA_{1,7}(q)\sms A_{2,6}(q) \big)
=q^{-1}\big( A_{7,1}(q)\sms qA_{6,2}(q) \big)\\
\aa_2(q)= q^{-2}\big(qA_{2,6}(q)\sms A_{3,5}(q) \big)
=q^{-2}\big( A_{6,2}(q)\sms qA_{5,3}(q) \big)  \\
\aa_3(q)= q^{-3}\big(qA_{3,5}(q)\sms A_{4,4}(q) \big)
= q^{-3}\big( A_{5,3}(q)\sms q A_{4,4}(q) \big)\\
\aa_4(q)=q^{-4} \big(qA_{4,4}(q)\big)
\end{align*}
\vskip-.2in
\noindent
From which we derive that
\begin{align*}
\aa_1(q)=   A_{1,7}(q)\sms A_{2,6}^{(1)}(q)
=  A_{7,1}^{(1)}(q)\sms  A_{6,2}(q)  \\
\aa_2(q)=   A_{2,6}^{(1)}(q)\sms A_{3,5}^{(2)}(q)
=  A_{6,2}^{(2)}(q)\sms  A_{5,3}^{(1)}(q)    \\
\aa_3(q)=  A_{3,5}^{(2)}(q)\sms A_{4,4}^{(3)}(q) \big)
=   A_{5,3}^{(3)}(q)\sms   A_{4,4}^{(2)}(q)  \\
\aa_4(q)=   A_{4,4}^{(3)}(q)\big)
\end{align*}
\vskip-.2in
\noindent
This gives
$$
   A_{4,4}=q^3\aa_4,\ess\ess
  \ess A_{3,5}= q^2(\aa_3+\aa_4),\ess\ess
  \ess A_{2,6}= q(\aa_2+\aa_3+\aa_4),\ess\ess
  \ess A_{1,7}= \aa_1+\aa_2+\aa_3+\aa_4,
$$
and
$$
A_{5,3}=q^3\aa_3+q^4\aa_4,\ess\ess\ess
A_{6,2}=q^5\aa_4+q^4\aa_3+q^2\aa_2,\ess\ess\ess
A_{7,1}=q^6\aa_4+q^5\aa_3+q^3\aa_2+q\aa_1
$$
Multiplying $A_{8-r,r}$ by $x^r$ and summing for $1\le r\le 7$
proves the odd case for $a=4$.
\sas

We will terminate by a statement with potentially 
significant  consequences.

\sas

\noindent{Conjecture IV}

{\ita The Polynomials $\aa_r(q)$ in 4.41 are unimodal and in $\BN[q]$.}

\end{document}